\documentclass[11pt,leqno]{article}
\usepackage{amsmath,amssymb,amsthm}

\numberwithin{equation}{section}

\newcommand{\cf}{{\it cf. \ }}
\newcommand{\ie}{{\it i.e. \ }}

\newcommand{\n}{\noindent}

\newcommand{\bb}[1]{\mathbb{#1}}
\newcommand{\cl}[1]{\mathcal{#1}}
\newcommand{\vp}{\varepsilon}

\theoremstyle{plain}
\newtheorem{thm}{Theorem}[section]
\newtheorem{defn}[thm]{Definition}
\newtheorem{cor}[thm]{Corollary}
\newtheorem{pro}[thm]{Proposition}
\newtheorem{klem}[thm]{Key Lemma}
\newtheorem{lem}[thm]{Lemma}

\theoremstyle{remark}
\newtheorem{rem}[thm]{Remark}
\newtheorem*{rk}{Remark}

\theoremstyle{definition}
\newtheorem*{nt}{Note}
\newtheorem*{note}{Notation}
\newtheorem{exm}[thm]{Example}

\begin{document}

\title{Are Unitarizable Groups Amenable?}

\author{by\\
Gilles Pisier\thanks{Supported in part by NSF   and
by the Texas Advanced Research Program 010366-163.}\\
Texas A\&M University\\
College Station, TX 77843, U. S. A.\\
and\\
Universit\'e Paris VI\\
Equipe d'Analyse, Case 186, 75252\\
Paris Cedex 05, France}

\date{Preliminary version, December 8,
2004}
\maketitle

\begin{abstract}
 We give a new
formulation
of some of our recent results on the following
problem: if all uniformly bounded representations
on a discrete group $G$ are similar to unitary
ones, is the group amenable? In   \S 5, we
give a new proof of Haagerup's theorem that, on
non-commutative  free groups, there are Herz-Schur
multipliers that are  not coefficients of
uniformly bounded representations. We actually
prove a refinement of this
result involving a generalization of the class
of Herz-Schur multipliers, namely the class
$M_d(G)$ which is formed of all the 
functions $f\colon \ G\to {\bb C}$ such that there are bounded
functions 
$\xi_i\colon \ G\to B(H_i, H_{i-1})$ ($H_i$ Hilbert) with $H_0 =
{\bb C}$, $H_d 
={\bb C}$ such that
$$f(t_1t_2\ldots t_d) = 
\xi_1(t_1) \xi_2(t_2)\ldots \xi_d(t_d).\qquad 
\forall~t_i\in G$$
We prove that if $G$ is a non-commutative  free
group, for any $d\ge 1$, we have
$$M_d(G)\not= M_{d+1}(G),$$ and hence there are
elements of
$M_d(G)$  which are  not coefficients
of uniformly bounded representations.
In the case $d=2$,  Haagerup's theorem implies
that $M_2(G)\not= M_{4}(G).$
\end{abstract}

\vfill\eject

{\bf Plan}

\S 0. Introduction

\S 1. Coefficients of uniformly bounded representations

\S 2. The spaces of multipliers $M_d(G)$

\S 3. The predual $X_d(G)$ of $M_d(G)$

\S 4. The $B(H)$-valued case

\S 5. A case study: The free groups.

\setcounter{section}{-1}

\section{Introduction}\label{sec0}

The starting point for this presentation is the following
result proved in the 
particular case $G={\bb Z}$ by Sz.-Nagy (1947).

\begin{thm}[Day, Dixmier 1950]\label{thm0.1}
Let $G$ be a locally
compact group. If 
$G$ is amenable, then every uniformly bounded (u.b.\ in short)
representation 
$\pi\colon \ G\to B(H)$ ($H$ Hilbert) is unitarizable. More
precisely, if we 
define
$$|\pi| = \sup \{\|\pi(t)\|_{B(H)}\mid t\in G\},$$
then, if $|\pi| < \infty$, there exists ${S}\colon\ H\to H$
invertible with 
$\|{S}\| \ \|{S}^{-1}\| \le |\pi|^2$ such that $t\to {S}^{-1}
\pi(t){S}$ is a 
unitary representation.
\end{thm}

\begin{nt}
 We say that  
$\pi\colon \ G\to B(H)$ ($H$ Hilbert)
 is unitarizable if there exists ${S}\colon\ H\to
H$ invertible such that $t\to {S}^{-1} \pi(t){S}$
is a  unitary representation.  We will mostly
restrict to discrete groups, but otherwise all 
representations $\pi\colon \ G\to B(H)$ are
implicitly assumed continuous on $G$  with
respect to the strong operator topology on
$B(H)$.
\end{nt}

\begin{defn}\label{defn0.2}
 We will say that a locally compact
group $G$ is 
unitarizable if every uniformly bounded (u.b.\ in short)
representation 
$\pi\colon \ G\to B(H)$ is unitarizable.
\end{defn}

In his 1950 paper, Dixmier \cite{Di} asked two questions which can be
rephrased as 
follows:
\begin{itemize}
\item[\bf Q1:] Is every $G$ unitarizable?
\item[\bf Q2:] If not, is it true that 
conversely unitarizable $\Rightarrow$ 
amenable?
\end{itemize}

\n In 1955, Ehrenpreis and Mautner answered Q1;
they showed that $G = {\rm SL}_2({\bb R})$  is not
unitarizable. Their work was clarified and
amplified in 1960 by Kunze-Stein \cite{KS1}.
See Remark \ref{rem0.7} below  for more recent work in this direction. 
 Here of
course $G = {\rm SL}_2({\bb R})$  is viewed as a Lie group,
but a fortiori the  discrete group $G_d$
underlying
   ${\rm SL}_2({\bb R})$ fails to be unitarizable, and
since every group is a quotient of a  free group
and ``unitarizable'' obviously passes to quotients, it
follows (implicitly)  that there is a
non-unitarizable free group, from which it is
easy to deduce  (since unitarizable passes to subgroups, see
Proposition \ref{pro0.5} below)
that
${\bb F}_2$ the free group with 2 generators is not
unitarizable. In the 80's,  many authors, notably
Mantero--Zappa \cite{MZ1}--\cite{MZ2}, Pytlik-Szwarc \cite{PyS}, 
Bo\.zejko--Fendler
\cite{BoF}, 
Bo\.zejko \cite{Bo2}, $\ldots$, and also M\l otkowski \cite{Ml}, Szwarc
\cite{Sz1}--\cite{Sz2}, Wysocza\'nski
\cite{W} (for free products of groups), produced explicit
constructions of u.b.\ non-unitarizable 
representations on ${\bb F}_2$ or on ${\bb F}_\infty$
(free group with countably  infinitely many
generators), see \cite{MPSZ} for a synthesis between
the Italian approach and the Polish one.
See  also Valette's papers \cite{V1}--\cite{V2} for the viewpoint of groups
acting on trees,
(combining Pimsner \cite{Pi} and \cite{PyS}) and \cite{Ja} for  recent work on
Coxeter groups.

\n This was partly
motivated by the potential  applications in
Harmonic Analysis of the resulting explicit
formulae (see e.g.\ \cite{FTP} and \cite{Ch}). For  instance, if we denote
by $|t|$ the length of an element in
${\bb F}_\infty$ (or in 
${\bb F}_2$), they constructed an analytic family $(\pi_z)_{z\in D}$
indexed by the 
unit disc in ${\bb C}$ of u.b.\ representations such that, denoting
by $\delta_e$ 
the unit basis vector of $\ell_2(G)$ at $e$ (unit element), we
have
\[
\forall~t\in G\qquad z^{|t|} = \langle\pi_z(t) \delta_e,
\delta_e\rangle.
\]
 Thus the
function $\varphi_z\colon \ t\to z^{|t|}$ is a
coefficient of a u.b.\  representation on
$G={\bb F}_\infty$. However, it can be shown that for
$z\notin 
{\bb R}$, the function $\varphi_z$ is {\it not\/} the coefficient of
any unitary 
representation, whence $\pi_z(\cdot)$ cannot be unitarizable.
A similar analysis can be made for the so-called spherical
functions.
 
Since unitarizable passes to subgroups (by``induction of
representations'', see
Proposition \ref{pro0.5} below) this implies

\begin{cor}\label{cor0.3}
 Any discrete group $G$ containing
${\bb F}_2$ as a subgroup 
is not unitarizable. 
\end{cor}

\begin{rem}\label{rem0.4}
 Thus if there is a discrete group $G$ which
is unitarizable but not 
amenable, this is a non-amenable group not containing ${\bb F}_2$.
The existence of 
such groups remained a fundamental open problem for many years
until Olshanskii
\cite{O1}--\cite{O2} established it in 1980, using  the
solution by Adian--Novikov (see \cite{A1}) of the
famous  Burnside problem, 
and also Grigorchuk's cogrowth criterion (\cite{Gri}). Later,
Adian (see
\cite{A2}) showed that the Burnside group 
$B(m,n)$ (defined as the universal group with $m$
generators such that every group element 
$x$ satisfies the relation $x^n=e$) are all
non-amenable when $m\ge 2$   and odd $n\ge 665$.
Obviously, since every element is periodic, such
groups cannot contain any  free infinite
subgroup. We should also mention   
 Gromov's    examples (\cite[\S 5.5]{Gr}) of 
  infinite discrete groups with Kazhdan's 
 property T (hence``very much" non-amenable) and
still without any
 free subgroup. In any case, it is natural to
wonder whether
the infinite Burnside groups are counterexamples
to the above Q2, whence the following.
\end{rem}

\n{\bf Question.} Are the Burnside groups $B(m,n)$
unitarizable?

In the next statement, we list the main stability 
properties of unitarizable groups.

\begin{pro}\label{pro0.5}
Let $G$ be a discrete group and let
$\Gamma$ be a  subgroup.
\begin{itemize}
\item[(i)] If $G$ is unitarizable, then $\Gamma$  also is.
\item[(ii)] If $\Gamma$ is normal, then $G$ is unitarizable only
if both
$\Gamma$  and $G/\Gamma$ are unitarizable. 
\end{itemize}
\end{pro}

\begin{proof}
\begin{itemize}
\item[(i)] Consider a u.b.\   representation $\pi\colon \ 
\Gamma\to B(H)$. By Mackey's induction, we have an ``induced''
representation 
$\hat\pi\colon \ G\to B(\widehat H)$ with $\widehat H\supset H$
that is  
 still u.b.\ (with the same bound) and hence is unitarizable.
Moreover, for any $t$ in $\Gamma$,  $\hat\pi(t)$ leaves $H$
invariant, and
 $\hat\pi(t)_{|H}=\pi(t)$. Hence, the original representation 
$\pi$ must also be unitarizable.
 (See \cite[p. 43]{P2} for full details). 
\item[(ii)] Let $q:\ G \to G/\Gamma$ be the quotient map and let
 $\pi\colon \ 
G/\Gamma\to B(H)$ be any representation. Then, trivially, $\pi$
is u.b.\  
(resp.\ unitarizable) iff  the same is true for $\pi \circ q$.
Hence
$G$ unitarizable implies $G/\Gamma$ unitarizable.
\end{itemize}
\end{proof}

\begin{rem}\label{rem0.6}
In (ii) above, we could not prove that
conversely
if $\Gamma$  and $G/\Gamma$ are unitarizable then $G$ is
(although the analogous
fact   for ideals in an operator algebra is true, see \cite[Exercise 27.1]{P7}).
In particular, we could not verify 
that the product of two unitarizable groups is 
unitarizable, however, it is known, and even for semi-direct
products,
if one of the groups is amenable, see \cite{NW}, see also \cite{P7} for
related
questions. 
Of course this should be true if unitarizable is the same as
amenable.
Similarly, it is not clear that a directed increasing union of a
family $(G_i)_{i\in I}$ of
 unitarizable groups
is unitarizable. Actually, we doubt that this is true in general.
However,
it is  true (and easy to check) if the family $(G_i)_{i\in I}$ is
``uniformly" unitarizable, in the following sense:
there is a function $F:\ {\bb R}_+ \to {\bb R}_+$  such that, for any
$i\in I$ and
any u.b.\ representation $\pi\colon \ G_i\to B(  H)$, there is
an invertible operator   ${S}\colon\ H\to H$   with 
$\|{S}\| \ \|{S}^{-1}\| \le F(|\pi|)$ such that $t\to {S}^{-1}
\pi(t){S}$ is a 
unitary representation.  
\end{rem}

\begin{rem}\label{rem0.7}
 There is an extensive literature continuing
  Kunze and Stein's work first  on ${\rm SL}(2,{\bb R})$ \cite{KS1} and
later
on ${\rm SL}(n,{\bb C})$ \cite{KS2}--\cite{KS4}, and devoted
(among other things) to the construction of
non-unitarizable uniformly bounded (continuous) 
representations on more general Lie groups.
We should mention P. Sally \cite{Sa1}--\cite{Sa2} for ${\rm SL}_2$ over local
fields
(see also \cite{MZ1}) and the universal covering group of
${\rm SL}(2,{\bb R})$, Lipsman \cite{Li1}--\cite{Li2} for the Lorentz groups
${\rm SO}_e(n,1)$ and for ${\rm SL}(2,{\bb C})$.
See the next remark for a synthesis of the current state of 
knowledge. 
We refer the reader to
 Cowling's papers (\cite{Co1,Co2}) 
for more recent work and
a much more comprehensive treatment of uniformly bounded
representations
on continuous groups. See also Lohou\'e's paper \cite{Lo}.
All in all, it seems there is a consensus among specialists
that discrete groups should be where to look primarily for a
counterexample
(i.e.\ unitarizable but not amenable),
if it exists. The next remark
hopefully should explain  why.
\end{rem}

\begin{rem}\label{rem0.72}
 (Communicated by Michael Cowling).
For an almost connected locally compact group~$G$ (that is,
$G/G_e$ 
is compact, where $G_e$ is the connected component of the
identity~$e$), 
unitarizability implies amenability.  The first step of the
argument 
for this is based on structure theory.  The group~$G$ has a
compact 
normal subgroup $N$ such that $G/N$ is a finite extension of a 
connected Lie group (see \cite[p.~175]{MZi}). Suppose that $G$ is
unitarizable. 
Then a fortiori $G/N$ is unitarizable.  If we can show that $G/N$
is 
amenable, then $G$ will be amenable, and we are done.  So we may
suppose 
that $G$~is a finite extension of a connected Lie group.  A
similar 
argument reduces to the case where $G$ is a connected Lie group,
and 
a third reduction (factoring out the maximal connected normal
amenable 
subgroup) leads to the case where $G$ is semisimple and
non-compact.  
It now suffices to show that a non-compact connected semisimple
Lie 
group~$G$ (which is certainly non-amenable) is not unitarizable.

So let $G$ be a non-compact connected semisimple Lie group.
We consider the representations $\pi_\lambda$ of~$G$ unitarily
induced 
from the characters $man \mapsto \exp(i \lambda \log a)$ of a
minimal
parabolic subgroup~$MAN$.  When $\lambda$ is real, $\pi_\lambda$
is 
unitary, and, according to B.~Kostant \cite{Ko}, $\pi_\lambda$ is 
unitarizable only if there is an element $w$ of the Weyl group 
$(\mathfrak{g},\mathfrak{a})$ such that $w \lambda =
\bar\lambda$.  
Take a simple root $\alpha$.  If $z$ is a complex number and
there 
exists $w$ in the Weyl group such that 
$w(z\alpha) = (z \alpha)\bar{\phantom{x}} = \bar z \alpha$, 
then $z$ is either purely real and $w\alpha = \alpha$ or $z$ is 
purely imaginary and $w\alpha = -\alpha$.  Thus if $z$ is neither 
real nor imaginary, then $\pi_{z\alpha}$ is not unitarizable.  
However, if the imaginary part of $z$ is small enough, then 
$\pi_{z\alpha}$ is uniformly bounded.  Indeed, using the
induction 
in stages construction (see \cite{Co2}, and also \cite{An}), we can make 
the representation uniformly bounded at the first stage, which 
involves a real rank one group only (see \cite{Co1} for the 
construction 
of the relevant Hilbert space) and then induce unitarily
thereafter 
to obtain a uniformly bounded representation. 
\end{rem}

The contents of this paper are as follows. In \S 1, we describe
our contribution 
on the above problem Q2, namely Theorem \ref{thm1.1} which says that if we
assume 
unitarizability with a specific quantitative bound then
amenability follows. We 
explain the main ideas of the proof in \S \ref{sec2}. There we introduce
our main objects 
of study in this paper namely the spaces $M_d(G)$.  The latter
are closely 
related on one hand to the space of ``multipliers of the Fourier
algebra,'' 
(which in our notation corresponds to $d=2$) and on the other
hand to the space 
$UB(G)$ of coefficients of {\it uniformly bounded} representations on
$G$, that we compare with the space $B(G)$ of
coefficients of {\it unitary} representations on
$G$. We have, for all $d\ge 2$
$$B(G)\subset UB(G) \subset M_d(G)\subset M_2(G).$$ Our
methods
lead  naturally to a new invariant of $G$, namely
the smallest $d$ such that $M_d(G) =  B(G)$,
that we denote by $d_1(G)$
(we set $d_1(G)=\infty$ if there is no such $d$). 
We have $d_1(G)=1$ iff $G$ is finite and
$d_1(G)=2$ iff $G$ is infinite and
amenable (see Theorem \ref{thm2.3}).
Moreover, we have  $d_1(G)=\infty$ 
when $G$ is any non Abelian free group.
Unfortunately, we cannot produce
any group with $2<d_1(G)<\infty$, and indeed such  an
example would provide a negative answer to the
above Q2.  While the main part of the paper is
partially expository, \S 5 contains  a
new result. We prove there that if
$G={\bb F}_\infty$ (free group with countably infinitely many
generators) then
$M_d(G)\ne M_{d+1}(G)$  for all $d\ge 2$. As a
corollary we obtain a completely different proof
of  Haagerup's unpublished result that $M_2(G)\ne
UB(G)$.

Let $H$ be a Hilbert space.
Actually, although it is less elementary, 
it is more natural to work with
the
$B(H)$-valued (or say $B(\ell_2)$-valued)
analogue of the spaces $M_d(G)$. The space $M_d(G)$
corresponds to $\dim(H)=1$ using ${\bb C}\simeq B({\bb C})$.
In the $B(\ell_2)$-valued case, the analogue
of $d_1(G)$ is denoted simply by $d(G)$.
We have  obviously $d_1(G)\le d(G)$ (but we do not
have examples where this inequality is strict).
Our results for $d(G)$ run parallel to those for $d_1(G)$.

 Although our methods (especially in the $B(\ell_2)$-valued case)
are inspired by the techniques of ``operator space theory'' and
``completely bounded  maps'' (see e.g.\ \cite{ER2},
\cite{Pa1} or
\cite{P8}]), we have strived to make our presentation 
accessible to a reader unfamiliar with those
techniques. This explains in particular why we present
the scalar valued (i.e. $\dim(H)=1$) case first.

We recall merely that a linear map $ u\colon\ A \to B(H)$
defined on a $C^*$-algebra $A$ is called completely bounded
(c.b.\ in short) 
if the maps $ u_n:\ M_n(A) \to M_n(B(H))$ defined by
$$ u_n( [a_{ij}]) = [u(a_{ij})] \quad  \forall [a_{ij}] \in
M_n(A)$$
are bounded uniformly over $n$ when $M_n(A)$
and $   M_n(B(H))$ 
 are each equipped with their unique  $C^*$-norm, i.e. 
the norm in the space of bounded operators acting on
$H\oplus\cdots\oplus H$ ($n$-times).

We also recall that, for any locally compact group $G$, 
 the  $C^*$-algebra of $G$ (sometimes called ``full" or
``maximal" 
to distinguish it from the ``reduced" case) is defined as the
completion
of the space $L_1(G)$ for the norm
defined by
$$\|f\| =\sup \left\| \int_G f(t) \pi(t) dt \right\|  \quad
\forall f\in L_1(G),$$
where the supremum runs over all 
(continuous) unitary representations $\pi$ on $G$.

In particular, the following result essentially due to Haagerup
(\cite{H2}) provides 
a useful (although somewhat  abstract) characterization of
unitarizable group 
representations.

\begin{thm}\label{thm0.8}
 Let $G$ be a locally compact group and let
$C^*(G)$ 
denote the (full) $C^*$-algebra of $G$. Let $\pi\colon \ G\to
B(H)$ be a 
uniformly bounded (continuous) representation. The following are
equivalent:
\begin{itemize}
\item[(i)] $\pi$ is unitarizable.
\item[(ii)] The mapping $ \tilde\pi\colon f\to \int f(t)\pi(t)dt$
from 
$L_1(G)$ to $B(H)$ extends to a completely bounded map from
$C^*(G)$ to $B(H)$.
\newline
  More generally, for an arbitrary bounded continuous function
$\varphi\colon \ G\to B(H)$,
the following are equivalent:
\item[(i)'] There is a unitary representation $\sigma\colon \
G\to B(H_\sigma)$ and
operators $\xi,\eta\colon\ H\to H_\sigma$ such that
$$\varphi(t)=\xi^* \sigma(t) \eta\qquad \forall t \in G.$$
\item[(ii)'] The mapping $\tilde\varphi\colon f\to \int
f(t)\varphi(t)dt$
extends to a completely bounded map from $C^*(G)$ to $B(H)$.
\end{itemize}
\end{thm}

\begin{proof}
 If $\pi$ is unitarizable, say we have $\pi(\cdot) = 
\xi\sigma(\cdot)\xi^{-1}$ with $\sigma$ unitary, then, by
definition of 
$C^*(G),\sigma$ extends to a $C^*$-algebra representation
$\hat\sigma$ from 
$C^*(G)$ to $B(H)$. Then, if we set $\hat\pi(\cdot) = 
\xi\hat\sigma(\cdot)\xi^{-1},\hat\pi$ extends $\pi$ and satisfies
(ii). Thus 
(i)~$\Rightarrow$~(ii). Conversely, if we have a completely
bounded extension 
$\hat\pi\colon \ C^*(G)\to B(H)$, then by \cite[Th. 1.10]{H2}, there
is $\xi$ invertible on 
$H$ such that $\xi^{-1}\hat\pi(\cdot)\xi$ is a $*$-homomorphism
(in other words 
a $C^*$-algebra morphism) and in particular
$t\to\xi^{-1}\hat\pi(t)\xi$ is a 
unitary representation, hence $\pi$ is unitarizable.
The proof of the equivalence of (i)' and (ii)' is analogous. That
(ii)' $\Rightarrow$ (i)' follows from the 
fundamental factorization of c.b. maps (see e.g.\ \cite[Chapt.
3]{P2}, \cite[p. 23]{P5},
or \cite{ER2}). The converse is obvious.
\end{proof}

In particular, this tells us that unitarizability is a countably
determined 
property:

\begin{cor}\label{cor0.9}
 Let $\pi\colon \ G\to B(H)$ be a
uniformly bounded 
representation on a discrete group $G$. Then $\pi$ is
unitarizable iff its
restriction to any countable subgroup $\Gamma\subset G$ is
unitarizable.
\end{cor}

\begin{proof}
If $\pi$ is not unitarizable, then, by Theorem \ref{thm0.8}, there is a
sequence $a^n$ 
with $a^n\in M_n(C^*(G))$ and $\|a^n\|\le 1$, $a^n_{ij}\in
\ell_1(G)$ and such 
that $\|[\tilde{\pi}(a^n_{ij})]\|_{M_n(B(H))}\to\infty$ when
$n\to\infty$. Since 
each entry $a^n_{ij}$ is countably supported, there is a
countable subgroup 
$\Gamma\subset G$ such that all the entries $\{a^n_{ij}\mid n\ge
1, 1\le i,j\le 
n\}$ are supported on $\Gamma$. This implies (by Theorem \ref{thm0.8}
again) that 
$\pi_{|\Gamma}$ is not unitarizable.  This proves the ``if''
part. The converse 
is trivial.
\end{proof}

\begin{cor}\label{cor0.10}
 If all  the countable subgroups of a
discrete group 
$G$ are unitarizable, then $G$ is unitarizable.
\end{cor}

\begin{proof}
This is an immediate consequence of the preceding corollary.
\end{proof}

\begin{rk}
 Let $G$ be a locally compact group and let
$\Gamma\subset G$ be 
a closed subgroup, we will say that $\Gamma$ is
$\sigma$-compactly generated if there is 
a countable union of compact subsets of $G$ that generates
$\Gamma$ as a closed 
subgroup. Then the preceding argument suitably modified shows
that, in the 
setting of Theorem~\ref{thm0.8}, if $\pi_{|\Gamma}$ is unitarizable for
any 
$\sigma$-compactly generated closed subgroup $\Gamma\subset G$,
then $\pi$ is 
unitarizable. 
\end{rk}

\section{Coefficients of uniformly bounded representations}\label{sec1}

It will be useful to introduce the space $B(G)$ of 
``coefficients of unitary representations''  
on a (discrete) group $G$ defined classically as 
follows.

We denote by $B(G)$ the space of all functions $f\colon \ G\to
{\bb C}$ for which 
there are a unitary representation
 $\pi\colon \ G\to B(H)$ and vectors $\xi,\eta\in
H$  such that
\begin{equation}\label{eq1.1}
\forall~t\in G\qquad f(t) = \langle \pi(t)\xi,\eta\rangle, 
\end{equation}
 This space can be equipped
with the norm
$$\|f\|_{B(G)} = \inf\{\|\xi\|\ \|\eta\|\}$$ 
where the infimum runs over all possible
$\pi,\xi,\eta$ as above. As is well  known, $B(G)$
is a Banach algebra for the pointwise product.
Moreover, $B(G)$  can be identified with the dual
of the ``full'' $C^*$-algebra of $G$, denoted by 
$C^*(G)$.

More generally, let $c\ge 1$ and let $G$ be a
semi-group with unit. In that case, we may replace
``representations' by unital semi-group homomorphisms.
Indeed, note that (since
a unitary operator is nothing but an invertible contraction
with contractive inverse)  a unitary representation on a group
$G$
is nothing but a unital semi-group homomorphism
$\pi\colon \ G\to B(H)$ such that 
$\sup \{\|\pi(t)\| \mid t\in G\}= 1$.
For any semi-group homomorphism
$\pi\colon \ G\to B(H)$, we again denote
$$|\pi|=\sup \{\|\pi(t)\| \mid t\in G\}.$$

In the sequel,  unless specified otherwise, $G$ will denote
a semi-group with unit.
  
 We denote by
$B_c(G)$ the space of all functions 
$f\colon \ G\to {\bb C}$ for which 
there is a  
unital semi-group homomorphism
$\pi\colon
\ G\to B(H)$ with $|\pi|\le c$ 
together with vectors
 $\xi,\eta$ in $H$ such that 
 \eqref{eq1.1} holds. 
Moreover, we denote
$$\|f\|_{B_c(G)} = \inf\{\|\xi\|\ \|\eta\|\mid f(\cdot) =
\langle\pi(\cdot)\xi, 
\eta\rangle \text{ with } |\pi|\le c\}.$$
Note that when $c=1$ and $G$ is a group, $B_1(G) =
B(G)$ with the same norm, since
$|\pi|=1$ iff $\pi$ is a unitary  
representation. For convenience of notation, we set
$B(G) =
B_1(G)$  also for semi-groups.
In the group case, $B_c(G)$ appears as
the space of coefficients of u.b.\ representations 
  with bound $\le c$.

The space $B_c(G)$ is a Banach space (for the above  norm).
Moreover, for any 
$c'\ge 1$ we have (since 
$\langle \pi(t)\xi,\eta\rangle \langle \pi'(t)\xi',\eta'\rangle
=\langle \pi\otimes \pi'(t)\xi\otimes \xi',\eta\otimes
\eta'\rangle$)
$$f\in B_c(G), g\in B_{c'}(G)\Rightarrow f\cdot g \in
B_{cc'}(G).$$
Note moreover that if $c\ge c'$ we have a norm one inclusion
$$B_{c'}(G) \subset B_c(G)$$
and in particular if $c'=1$ we find
$$\|f\|_{B_c(G)} \le \|f\|_{B(G)}.\leqno \forall~f\in B(G)$$

We will denote by $ {UB}(G)$ the space of coefficients of 
uniformly bounded representations on $G$;
in other words, we set:
$$ {UB}(G)=\bigcup_{c>1}
B_c(G).$$

The following result partially answers Dixmier's question Q2.

\begin{thm}[\cite{P3}]\label{thm1.1}
The following
properties of a discrete group $G$ are 
equivalent.
\begin{itemize}
\item[(i)] $G$ is amenable.
\item[(ii)] $\exists K \ \exists \alpha<3$ such that for every
u.b.\ 
representation $\pi\colon \ G\to B(H)$ $\exists{S}\colon \ H\to
H$ invertible 
with $\|{S}\|\ \|{S}^{-1}\| \le K|\pi|^\alpha$ such that
${S}^{-1}\pi(\cdot){S}$ 
is a unitary representation.
\item[(ii)$'$] Same as (ii) with $\alpha=2$ and $K=1$.
\item[(iii)] $\exists K\ \exists \alpha<3$ such that for any
$c>1$ $B_c(G) 
\subset B(G)$ and we have
\[
\forall~f\in B_c(G)\qquad
\|f\|_{B(G)} \le Kc^\alpha \|f\|_{B_c(G)}.
\]
\item[(iii)$'$] Same as (iii) with $\alpha=2$ and $K=1$.\medskip
\end{itemize}
\end{thm} 

\begin{rk}
 Actually, the preceding result
remains valid for a general locally compact group.
Indeed, as Z.J.\ Ruan told us, the 
Bo\. zejko criterion \cite{Bo1} which we use to prove
Theorem \ref{thm1.1} for discrete groups 
(which says, with the notation explained below,
that $M_2(G)=B(G)$ iff $G$ is amenable)
remains valid in the general case. Z.J.\ Ruan
checked that Losert's proof of a similar but a
priori weaker statement specifically for the
non-discrete case (see \cite{Lo}) can be modified to
yield this, and apparently (cf.\ \cite{Ru}) this fact
was already known to Losert (unpublished). Our
proof that (iii) above implies
 $M_2(G)=B(G)$ does not really use the
discreteness
of the group, whence the result in full
generality. 
\end{rk}

\n This observation, concerning the extension
to
 general locally compact groups, was also made independently
by Nico Spronk \cite{Spr}.

\n We would like to emphasize that there are two separate,
very different arguments: one for the discrete
case and one for the    non-discrete one. This
is slightly surprising. A unified approach would
be interesting.

\begin{proof}[First part of the proof of Theorem \ref{thm1.1}]
That  (i) $ \Rightarrow$  (ii)$'$ is the Dixmier-Day
 result mentioned at the beginning.
The implications (ii)$'
\Rightarrow$ (ii) and (iii)$' \Rightarrow$ (iii)
are  trivial. Moreover, (ii) $\Rightarrow$ (iii)
and (ii)$'
\Rightarrow$ (iii)$'$ are  easy to check:\ indeed
consider $f(\cdot) = \langle \pi(\cdot)
\xi,\eta\rangle$  with $|\pi|\le c$. Then, if $S$
is such that
$\hat \pi = S\pi(\cdot)S^{-1}$ is a unitary 
representation, we have $f(\cdot) =
\langle\pi(\cdot)\xi,\eta\rangle = 
\langle\hat\pi(\cdot) S\xi, (S^{-1})^*\eta\rangle$, 
hence the coefficients of 
$\pi$ are coefficients of $\hat\pi$ and 
$$\|f\|_{B(G)}\le \|S\xi\|\ 
\|(S^{-1})^*\eta\| \le\|S\|\ \|S^{-1}\|\ \|\xi\|\
\|\eta\|$$
 whence
$$\|f\|_{B(G)} \le \|S\|\ \|S^{-1}\| \ \|f\|_{B_c(G)}.$$
Now if (ii) holds we can find $S$ as above with $\|S\|\
\|S^{-1}\|\le 
Kc^\alpha$, thus  (ii) $\Rightarrow$ (iii), and similarly (ii)$'
\Rightarrow$ 
(iii)$'$.
Thus it only remains to prove
that   (iii) $ \Rightarrow$  (i).
\end{proof}

\section{The spaces of multipliers $\pmb{M_d(G)}$}\label{sec2}

To explain the proof of Theorem \ref{thm1.1}, we will need
some additional notation.

\begin{note}
 Let $d\ge 1$ be an integer.
Let $G$ be a a semigroup with unit.
We are mainly interested in the group case, but we could 
also take $G={\bb N}$. 
\end{note} 

\n  Let
$M_d(G)$ be
 the space of all 
functions $f\colon \ G\to {\bb C}$ such that there are bounded
functions 
$\xi_i\colon \ G\to B(H_i, H_{i-1})$ ($H_i$ Hilbert) with $H_0 =
{\bb C}$, $H_d 
={\bb C}$ such that
\begin{equation}\label{eq2.1}
\forall~t_i\in G\qquad f(t_1t_2\ldots t_d) = \xi_1(t_1) \xi_2(t_2)\ldots
\xi_d(t_d).
\end{equation}
Here of course we use the identification $B(H_0,H_d) =
B({\bb C},{\bb C})\simeq {\bb C}$. We 
define
$$\|f\|_{M_d(G)} = \inf\{\sup_{t_1\in G} \|\xi_1(t_1)\|\ldots
\sup_{t_d\in G} 
\|\xi_d(t_d)\|\}$$
where the infimum runs over all possible ways to write $f$ as in 
\eqref{eq2.1}.

The definition of the spaces ${M_d(G)}$ and of
the  more general spaces ${M_d(G;H)}$ appearing
below is motivated by the work of
Christensen-Sinclair on ``completely
bounded multilinear maps" and the so-called
Haagerup tensor product (see \cite{ChS}). The
connection is explained in detail in \cite{P3},
and is important for the proofs
of all the results below, but we prefer
to skip this in the present exposition (see
however \S \ref{sec4} below). 

When $d=2$, and $G$ is a group, the space $M_2(G)$
is the classical space of ``Herz-Schur 
multipliers'' on
$G$. This space also coincides (see  \cite{BoF} or
\cite[p. 110]{P2}) with the  space of all c.b.\
``Fourier multipliers'' on the reduced
$C^*$-algebra 
$C^*_\lambda(G)$. The question whether 
the space $M_2(G)$ coincides with the 
space of coefficients of u.b.\ representations
 (namely $\bigcup\limits_{c>1} 
B_c(G)$) remained open for a while but Haagerup \cite{H1} showed that
it is not the 
case. More precisely, he showed that if $G = {\bb F}_\infty$, we have
\[
\forall~c>1\qquad B_c(G) \underset{\ne}{\subset}
M_2(G).
\]
We give a different proof of a more precise
statement in \S \ref{sec5} below.

For $d>2$, in the group case, 
the spaces $M_d(G)$ are not so naturally interpreted 
in terms of ``Fourier" multipliers. In particular, in spite of
the strong analogy with the multilinear multipliers
introduced in \cite{ER1} (those are complex valued functions
on $G^d$), there does not seem to be any significant connection.

In the case $G={\bb N}$, the space $M_3(G)$ is characterized
in \cite{P6} as the space of ``completely shift bounded"
Fourier multipliers on
the Hardy space $H_1$, but this interpretation
is restricted to $d=3$ and uses the commutativity.

\begin{rem}\label{rem2.1}
 Note the following easily
checked
  inclusions,  valid when $G$ is a group or a
semigroup with unit:
\begin{align*}
B(G)=B_1(G) &\subset UB(G)=\bigcup_{c>1}
B_c(G)
\subset M_d(G) \subset  M_{d-1}(G) 
\subset\cdots\\
\cdots &\subset M_2(G) \subset M_1(G) =
\ell_\infty(G),
\end{align*}
 and we have clearly
\begin{equation}\label{eq2.2}
\forall m\le d\qquad \|f\|_{M_{m}(G)} \le
\|f\|_{M_d(G)}.
\end{equation}
Moreover, we have
\begin{equation}\label{eq2.3}
\forall~f\in B_c(G)\qquad \|f\|_{M_d(G)} \le c^d\|f\|_{B_c(G)}.
\end{equation}
  Indeed, if
$f(\cdot) = \langle\pi(\cdot)\xi,\eta\rangle$ with
$|\pi|\le c$, then  we can write
\begin{align*}
f(t_1t_2\ldots t_d) &= \langle \pi(t_1)\ldots 
\pi(t_d)\xi,\eta\rangle\\
&= \xi_1(t_1)\xi_2(t_2)\ldots \xi_d(t_d)
\end{align*}
where $\xi_1(t_1) \in B(H_\pi,{\bb C})$, $\xi_d(t_d) \in
B({\bb C},H_\pi)$ and 
$\xi_i(t_i) \in B(H_\pi,H_\pi)$ $(1<i<d)$ are defined by
$\xi_1(t_1)h = 
\langle\pi_1(t_1)h,\eta\rangle$ $(h\in H_\pi)$ $\xi_d(t_d)\lambda
= 
\lambda\pi(t_d)\xi$ $(\lambda\in{\bb C})$ and $\xi_i(t_i)= \pi(t_i)$
$(1<i<d)$. 
Therefore, we have
\begin{align*}
\|f\|_{M_d(G)} &\le \sup\|\xi_1\|
\sup\|\xi_2\|\ldots\sup 
\|\xi_d\|\\
&\le |\pi|^d\|\xi\|\ \|\eta\|\le c^d\|\xi\|\ \|\eta\|
\end{align*}
whence the announced inequality  \eqref{eq2.3}.
\end{rem}

\begin{rem}\label{rem2.2}
 It is easy to see 
(using tensor products) that $M_d(G)$
is a unital Banach algebra for the pointwise product of
functions on $G$: for any $f,g$ in $M_d(G)$ we have
$$ \|fg\|_{M_d(G)} \le \|f\|_{M_d(G)}\|g\|_{M_d(G)}.$$
The function identically equal to 1 on $G$ is the unit and has norm 1.
\end{rem}

\begin{rk}
Let $h\colon \ \Gamma\to G$ be a unital
homomorphism between 
two  groups (or two semi-groups with unit). Then for any $f$ in
$M_d(G)$ the 
composition $f\circ h$ is in $M_d(\Gamma)$ with $\|f\circ
h\|_{M_d(\Gamma)} \le 
\|f\|_{M_d(G)}$. The proof is obvious.
\end{rk}

In particular, if $\Gamma\subset G$ is a subgroup we have 
$\|f_{|\Gamma}\|_{M_d(\Gamma)}\le \|f\|_{M_d(G)}$. Moreover, if
$\Gamma$ is a 
normal subgroup in a group $G$ and if $q\colon \ G\to G/\Gamma$
is the quotient 
map, then $\|f\|_{M_d(G/\Gamma)} = \|f\circ q\|_{M_d(G)}$.
(Indeed, the equality 
can be proved easily using an arbitrary pointwise lifting
$\rho\colon \ G/\Gamma 
\to G$.) Let $\Gamma\subset G$ be again an arbitrary subgroup of
a group $G$. 
Given a function $f\colon \ \Gamma\to {\bb C}$, we let $\tilde
f\colon \ G\to {\bb C}$ 
be the extension of $f$ vanishing outside $\Gamma$. Then, it is
rather easy to 
see that $\|\tilde f\|_{M_2(G)} = \|f\|_{M_2(\Gamma)}$ (but the
analogue of this 
for $d>2$ seems unclear). 

\n It is well known that $1_\Gamma$ is in the unit ball 
of $B(G)$ (hence a fortiori of $M_d(G)$) hence by Remark~\ref{rem2.2} we
have for any $f$ 
in $M_d(G)$:
$$\|f\cdot 1_\Gamma\|_{M_d(G)} \le \|f\|_{M_d(G)}.$$

The proof of (iii) $\Rightarrow$ (i) in Theorem \ref{thm1.1}
 uses the following criterion
 for amenability due to Marek Bo\.zejko 
\cite{Bo1}.

\begin{thm}\label{thm2.3}
 Let $G$ be a discrete group.
Then $G$ is amenable iff $B(G)  = M_2(G)$.
\end{thm}

\begin{rk}
We do not know whether $B(G) = M_3(G)\Rightarrow
G$ amenable. 
\end{rk}

\begin{proof}[Sketch of proof of Theorem \ref{thm2.3}]
 The only if
part is quite easy. Let us  sketch the proof of
the ``if'' part. Assume $B(G) = M_2(G)$. Then
there is a  constant $K$ such that, for any $f$
in the space ${\bb C}[G]$ of all finitely  supported
functions $f\colon \ G\to {\bb C}$, we have
$$\|f\|_{B(G)} \le K\|f\|_{M_2(G)}.$$
Let $\varepsilon\colon \ G\to \{-1,1\}$ be a ``random choice of
signs'' indexed 
by $G$, and let ${\bb E}$ denote the expectation with respect to the
corresponding 
probability. We will estimate the average of the norms of the
pointwise product 
$\varepsilon f$. More precisely we claim that there are numerical
constants $C'$ 
and $C''$ (independent of $f$) such that
\begin{align}
\label{eq2.4}
\left(\sum_{t\in G} |f(t)|^2\right)^{1/2} &\le
C'{\bb E}\|\varepsilon 
f\|_{B(G)}\\
\label{eq2.5}
{\bb E}\|\varepsilon f\|_{M_2(G)} &\le C''\left\|\sum |f(t)|^2 
\lambda(t)\right\|^{1/2}_{C^*_\lambda}.
\end{align}
Using this it is easy to conclude:\ indeed we have
$$\sum|f(t)|^2 \le (C' K C'')^2 \left\|\sum |f(t)|^2 
\lambda(t)\right\|_{C^*_\lambda}$$
and by the well known Kesten-Hulanicki criterion
(cf.\ e.g.\ \cite[Th. 2.4]{P2}), this implies that $G$ is  amenable.

We now return to the above claims. The inequality
\eqref{eq2.4} can be seen as a  consequence of the fact (due
to N.~Tomczak-Jaegermann \cite{TJ}) that $B(G)$ is of
cotype  2 (a Banach space $B$ is called of cotype 2
if there is a constant $C$ such that for any finite sequence
$(x_i)$ in $B$, the following inequality holds
$(\sum \|x_i\|^2 )^{1/2} \le C {\rm Average}_{\pm} \|\sum \pm 
x_i\|$). 

 As for  \eqref{eq2.5}, it is proved
in \cite{Bo1} using an idea due to Varopoulos \cite{Va}.  However, more
recently the following result was proved in \cite{P1}:\
Consider all  possible ways to have the following
decomposition
\begin{equation}\label{eq2.6}
   f(t_1t_2) = \alpha(t_1,t_2) + \beta(t_1,t_2)
\qquad \forall~t_1,t_2\in G \end{equation}
and let 
$$|||f||| = \inf\left\{\sup_{t_1}
\left(\sum_{t_2}  |\alpha(t_1,t_2)|^2\right)^{1/2}  +
\sup_{t_2} \left(\sum_{t_1} |\beta(t_1,t_2)|^2 
\right)^{1/2}\right\}$$
where the infimum runs over all possible decompositions as in 
\eqref{eq2.6}.

\n Then (cf.\ \cite{P1}) there is a numerical constant
$\delta>0$ such that 
\begin{equation}\label{eq27}
\delta|||f||| \le {\bb E}\|\varepsilon 
f\|_{M_2(G)} 
\le
|||f|||\qquad\forall~f\in {\bb C}[G] .
\end{equation}

Note that the right-hand side of \eqref{eq27} is an immediate consequence of 
the following inequality
\begin{equation}\label{eq2.88}
\|f\|_{M_2(G)} \le |||f|||,
\end{equation}
and the latter is easy:\ we simply write
\[
f(t_1t_2) = \langle\xi_1(t_1), \xi_2(t_2)\rangle + \langle \eta_1(t_1), 
\eta_2(t_2)\rangle
\]
where 
\begin{align*}
\xi_1(t_1) &= \sum\limits_{t_2} \alpha(t_1,t_2) \delta_{t_2}, \quad \xi_2(t_2)  = 
\delta_{t_2}\\
\intertext{and}
\eta_2(t_2) &= \sum_{t_1} \overline {\beta(t_1,t_2)} \delta_{t_1},\quad \eta_1(t_1) = 
\delta_{t_1},
\end{align*}
and \eqref{eq2.88} follows.

Moreover, we also have
\begin{equation}\label{eq2.89}
|||f||| \le \left\|\sum |f(t)|^2 \lambda(t)\right\|^{1/2}_{C^*_\lambda},
\end{equation}
therefore \eqref{eq2.5} follows from \eqref{eq27} and \eqref{eq2.89} with $C''=1$. The 
inequality \eqref{eq2.89} follows from the following observation:\ $\|\sum 
|f(t)|^2 \lambda(t)\|_{C^*_\lambda} \le 1$ iff there is a decomposition of the 
form
\begin{equation}\label{eq2.90}
|f(t_1t_2)| = |a(t_1,t_2)|^{1/2} \cdot |b(t_1,t_2)|^{1/2}\qquad \forall 
t_1,t_2\in G
\end{equation}
for kernels $a,b$ on $G\times G$ such that
\begin{equation}\label{eq2.91}
\sup_{t_1} \sum_{t_2} |a(t_1,t_2)|^2 \le 1\quad \text{and}\quad \sup_{t_2} 
\sum_{t_1} |b(t_1,t_2)|^2 \le 1.
\end{equation}
 Then \eqref{eq2.90} and \eqref{eq2.91} imply that $|||f|||\le 1$, and hence 
\eqref{eq2.89} follows by homogeneity. Indeed, by a compactness argument,
these assertions are
immediate consequences of the following  Lemma.

\end{proof}

This Lemma gives a converse to Schur's classical criterion for boundedness
on $\ell_2$ of matrices with positive entries
 (we include the proof for lack of a suitable reference).
 See \cite{P55} for more information on this.
 \begin{lem}\label{lem2.44}
Let $n\ge 1$. Let $\{f_{ij}\mid 1\le i,j\le n\}$ be complex scalars such that 
the matrix 
$[|f_{ij}|^2]$ has norm $\le 1$ as an operator on the Euclidean space $\ell_2^n$.
Then there are $(a_{ij})$ and $(b_{ij})$ with
\[
\sup_i \sum\nolimits_j|a_{ij}|^2 \le 1\quad \text{and}\quad \sup_j \sum\limits_i 
|b_{ij}|^2 \le 1
\]
such that $|f_{ij}|  = |a_{ij}|^{1/2} |b_{ij}|^{1/2}$. Therefore, there are 
$(\alpha_{ij})$ and $(\beta_{ij})$ with
\[
\sup_i (\sum\nolimits_j |\alpha_{ij}|^2)^{1/2} + \sup_j (\sum\nolimits_i |\beta_{ij}|^2)^{1/2} \le 1
\]
such that $f_{ij} = \alpha_{ij} + \beta_{ij}$.
\end{lem}

\begin{proof}
By perturbation and compactness arguments, we can assume that $|f_{ij}|>0$ for 
all $i,j$. Let $T = [|f_{ij}|^2]$. We may assume $\|T\| =1$. Let $\xi = (\xi_i)$ 
be a Perron--Frobenius vector for $T^*T$ so that $\xi_i>0$ for all $i$ and 
$T^*T\xi = \xi$. Let $\eta = T\xi$, so that $T^*\eta=\xi$. If we then set 
$|a_{ij}|^2 = |f_{ij}|^2 \xi_j\eta^{-1}_i$ and $|b_{ij}|^2 = |f_{ij}|^2\xi^{-1}_j\eta_i$ we 
obtain the first assertion. 
By the arithmetic-geometric mean inequality we have
 $|f_{ij}|\le g_{ij}$ with $g_{ij}=2^{-1}(|a_{ij}|+|b_{ij}|)$. If we then set
\[ \alpha_{ij}=2^{-1}|a_{ij}|f_{ij}g_{ij}^{-1} 
 \quad \text{and}\quad \beta_{ij}=2^{-1}|b_{ij}|f_{ij}g_{ij}^{-1} ,
\]
we obtain the second assertion.
\end{proof}

\begin{rk}
 The proof of Theorem \ref{thm2.3} sketched
above shows that $G$ is  amenable if $\exists K$
$\forall f \in {\bb C}[G]$
$$\|f\|_{B(G)} \le K\|f\|_{M_2(G)}.$$
\end{rk}

\begin{rk}
   Note that \eqref{eq27}
and  \eqref{eq2.88} show that
\[\sup_{\varepsilon} \|\varepsilon 
f\|_{M_2(G)}  \le \delta^{-1}  {\bb E}\|\varepsilon 
f\|_{M_2(G)}. 
\]
\end{rk}

The proof of the implication (iii) $\Rightarrow$ (i) in
Theorem~\ref{thm1.1} rests on the 
following.

\begin{klem}[Implicit in \cite{P3}]\label{lem2.4}
Let
$f\in B(G)$. Fix $d\ge 1$. Then, for any 
$c\ge 2$, we have
$$\|f\|_{B_c(G)} \le 2 \|f\|_{M_d(G)} + 2c^{-(d+1)}
\|f\|_{B(G)}.$$
More generally, for any $1\le \theta <c$, we have for any $f$ in
$B_\theta(G)$ 
and any $d\ge 1$
$$ {\|f||_{B_c(G)}  \le \left(\sum^d_{m=0}
(\theta/c)^m\right) \cdot 
\|f\|_{M_m(G)} 
 \quad + \left(\sum_{m>d}
(\theta/c)^m\right)\cdot \|f\|_{B_\theta(G)}.}$$
\end{klem}

\begin{rk}
 The proof of the key lemma uses ideas from two
remarkable 
papers due to Peller \cite{Pe} and Blecher and Paulsen \cite{BP}.
\end{rk}

\begin{proof}[Proof of (iii) $\Rightarrow$ (i) in Theorem \ref{thm1.1}]
 Assume
(iii). Then using 
the key lemma with $d=2$ (and $\theta=1$), we have
for all
$f$ in
$B(G)$ and all $c\ge 2$
\begin{align*}
\|f\|_{B(G)} &\le Kc^\alpha \|f\|_{B_c(G)}\\
&\le 2Kc^\alpha \|f\|_{M_2(G)} + 2Kc^{\alpha-3} \|f\|_{B(G)}.
\end{align*}
But we can choose $c = c(K,\alpha)$ large enough
so that $2Kc^{\alpha-3} = 1/2$  (say) and then we
obtain
$$\left(1-\frac12\right) \|f\|_{B(G)} \le 2Kc^\alpha
\|f\|_{M_2(G)}$$
so that we conclude
$$\|f\|_{B(G)} \le 4K c(K,\alpha)^\alpha \|f\|_{M_2(G)},$$
hence, by  Theorem \ref{thm2.3}, $G$ is amenable. 
\end{proof}

The proof of the key lemma is based on the following result (of
independent 
interest) which is ``almost'' a characterization of $B_c(G)$.

\begin{thm}\label{thm2.5}
 Fix a number $c\ge 1$.
Consider $f\in \bigcap\limits_{m\ge  1} M_m(G)$
such that   $\sum_m c^{-m}
\|f\|_{M_m(G)}  < \infty$. Then $f\in
B_c(G)$  and moreover
\begin{equation}\label{eq2.8}
\|f\|_{B_c(G)} \le |f(e)| + \sum_{m\ge 1} c^{-m}
\|f\|_{M_m(G)}.
\end{equation}
Conversely, for all $f$ in $B_c(G)$, we have
\begin{equation}\label{eq2.9}
\sup_{m\ge 1} c^{-m} \|f\|_{M_m(G)} \le
\|f\|_{B_c(G)}.
\end{equation}
\end{thm}

\begin{nt}
\eqref{eq2.9} is easy and has been proved already (see 
\eqref{eq2.2}). The main 
point is  \eqref{eq2.8}. 
\end{nt}

\begin{proof}
 This is essentially \cite[Theorem 1.12]{P3} and the remark following it. For 
the
convenience of the reader, we give some more
details.

\n  In \cite{P3} the natural predual of ${B_c(G)}$ is
considered and denoted by $\widetilde A_c$. By
\cite[Th. 1.7]{P3}, any $x$ in the open unit ball
of $\widetilde A_c$ can be written
as
$$x=\sum_{m=0}^\infty c^{-m} x_m$$
where each $x_m$ is an element of $C^*(G)$ 
which is the image, under the natural product
map
of an element $X_m$ in the   unit ball
of $\ell_1(G) \otimes_h\cdots \otimes_h
\ell_1(G)$ ($m$ times). This implies
by duality,  for all $m> 0$ 
$$  |\langle f,x_m\rangle |\le \|f\|_{M_m(G)}.$$
In the particular case $m=0$,  $x_0$ is
a multiple of the unit $\delta_e$ by a scalar
of modulus $\le 1$.
Whence
$$ |\langle
f,x\rangle|\le |f(e)|+ \sum_{m=1}^\infty
c^{-m}\|f\|_{M_m(G)}$$
and a fortiori, by  \eqref{eq2.2} and  \eqref{eq2.3} 
$$\le |f(e)|+\sum_{m=1}^d c^{-m} \|f\|_{M_d(G)}
+ \sum_{m=d+1}^\infty c^{-m} \|f\|_{B(G)}.$$
\end{proof}
 
\begin{cor}\label{cor2.6}
 Let $UB(G) =
\bigcup\limits_{c>1} B_c(G)$ be the  space of
coefficients of u.b.\ representations on $G$.
Then $f\in UB(G)$ iff 
$\sup\limits_{m\ge 1} \|f\|^{1/m}_{M_m(G)} <
\infty$. More precisely, let $c(f)$ denote the
infimum of the numbers $c\ge 1$ for which
$f\in B_c(G)$. Then, we have
$$c(f)=\limsup\limits_{m\to \infty}
\|f\|^{1/m}_{M_m(G)}.$$
\end{cor}

\begin{proof}[Proof of Key Lemma \ref{lem2.4}]
 This is an easy
consequence of  \eqref{eq2.8}, \eqref{eq2.9} and the  obvious
inequalities
\begin{align*}
|f(e)| &\le \|f\|_{\ell_\infty(G)} \le \|f\|_{M_2(G)}
\le\cdots\le 
\|f\|_{M_d(G)}\\
&\le\cdots\le \|f\|_{B(G)}.\qquad\qed
\end{align*}
\renewcommand{\qed}{}\end{proof}

In the case $c=1$, Theorem \ref{thm2.5} seems to degenerate but actually
the following 
``limiting case'' can be established, as a rather simple
dualization of a result 
in \cite{BP}.

\begin{pro}\label{pro2.7}
 Consider a function
$f\colon \ G\to {\bb C}$. Then $f\in  B(G)$ iff $f\in
\bigcap\limits_{m\ge 1} M_m(G)$ with
$\sup\limits_m 
\|f\|_{M_m(G)} < \infty$. Moreover we have
$$\|f\|_{B(G)} = \sup_{m\ge 1} \|f\|_{M_m(G)}.$$
\end{pro}

\begin{rk}
The same argument shows the following. Given a
real number 
$\alpha\ge 0$, we say that $G$ satisfies the condition
$(C_\alpha)$ if there is 
$K\ge 0$ such that for any $f$ in $B(G)$ we have
\[
\forall~c>1\qquad 
\|f\|_{B(G)} \le Kc^\alpha\|f\|_{B_c(G)}.
\]
Then the preceding argument shows that if $d\le \alpha < d+1$,
$(C_\alpha)$ 
implies that $B(G) = M_d(G)$ (with equivalent norms).
\end{rk}

We are thus led to define the following
quantities:
\begin{align*}
\alpha_1(G) &= \inf\{\alpha\ge 0\mid G \text{
satisfies } 
(C_\alpha)\}\\
d_1(G) &= \inf\{d\in {\bb N} \mid M_d(G) = B(G)\}.
\end{align*}
With this notation, the preceding argument shows
that $d_1(G)\le \alpha_1(G)$.
A priori $\alpha_1(G)$ is a real number, but
(although we have no direct  argument for this)
it turns out that it is an integer:

\begin{thm}\label{thm2.8}
 Assume $B(G) =
\bigcup\limits_{c>1} B_c(G)$. Then 
$\alpha_1(G) <\infty$ and moreover
$$\alpha_1(G) = d_1(G).$$
In particular, we have 
$$  M_d(G)=
M_{d+1}(G) \quad \forall d\ge \alpha_1(G) .$$
\end{thm}

Actually, for the last assertion to hold, it
suffices to have much less:

\begin{thm}[\cite{P4}]\label{thm2.9}
Let $G$ be a
semigroup with unit.  Assume that there are
$1\le \theta <c$ such that $B_{\theta}(G) =
B_c(G)$. Then there is an integer $D$ such that
$B_{\theta}(G) =M_D(G)$, and in particular, we
have 
$$ M_d(G)=
M_{d+1}(G) \quad \forall d\ge D.$$
\end{thm}

\begin{rk}
Let $G$ be a locally compact group.
Let $G_d$ be $G$ equipped with the discrete topology. 
In \cite{H3}, Haagerup proves that if a function
 $\phi\colon\ G\to {\bb C}$ belongs to $M_2(G_d)$ and
is continuous, then it belongs to $M_2(G)$
(with the same norm). We do not know if the analogous 
statement is valid for $M_3(G)$ or $M_d(G)$ 
when
$d\ge3$. 
\end{rk}

\begin{rem}\label{rem2.10}
Let $I_1,\ldots, I_d$ be arbitrary sets.
We will denote 
by 
$M_d(I_1,\ldots, I_d)$ the space of all functions $f\colon \ 
I_1\times\cdots\times I_d\to {\bb C}$ for which there are bounded
functions 
$f_i$
$$f_i\colon \ I_i \to B(H_i,H_{i-1})\quad
(\text{here}\  H_i \ \text{are Hilbert spaces}\ \text{with } H_d =
H_0 = {\bb C}) \text{  such
that}$$
$$\forall~b_i\in I_i\qquad f(b_1,\ldots, b_d) = f_1(b_1)\ldots f_d(b_d).$$
We equip this space with the norm
$$\|f\| = \inf\left\{\prod^d_{i=1}
 \sup_{b\in I_i} \|f_i(b)\|\right\}$$
where the infimum runs over all possible such factorizations.
\end{rem}

 In particular, if $I_1 = I_2 = \cdots = I_d = G$, then for 
any function $\varphi$ in $M_d(G)$, we have
$$\|\varphi\|_{M_d(G)} = \|\Phi\|_{M_d(G,\ldots, G)}$$ 
where $\Phi$ is defined by
$$\Phi(t_1,\ldots, t_d) = \varphi(t_1t_2\ldots t_d).$$
By a well known trick, one can check that $f\to
\|f\|_{M_d(I_1,\ldots, I_d)}$ is 
subadditive (and hence is a norm) on $M_d(I_1,\ldots, I_d)$, i.e.
that  we have for all $f,g$ in   $M_d(I_1,\ldots, I_d)$ :
\begin{equation}\label{subadd}
\| f+g\|_{M_d(I_1,\ldots, I_d)}\le \| f \|_{M_d(I_1,\ldots, I_d)}+\|
 g\|_{M_d(I_1,\ldots, I_d)}
\end{equation}
 Let
us quickly 
sketch this:\ Let $f,g$ be in the open unit ball of
$M_d(I_1,\ldots, I_d)$. Then 
by homogeneity we can assume
\begin{equation}\label{sub}
f(b_1,\ldots, b_d) = f_1(b_1)\ldots f_d(b_d)\quad \text{and}\quad 
g(b_1,\ldots, b_d) = g_1(b_1)\ldots g_d(b_d)
\end{equation}
with $\sup\|f_j(b)\| < 1$ and $\sup\|g_j(b)\| < 1$ for all $j$.
Then we can 
write for any $0\le\alpha \le 1$
$$(\alpha f + (1-\alpha)g) (b_1,\ldots, b_d) = F_1(b_1)\ldots
F_d(b_d)$$
where
\begin{gather*}
F_1(b_1) = [\alpha^{1/2}f_1(b_1) \quad (1-\alpha)^{1/2}
g_1(b_1)]
\quad \text{(row matrix with operator entries)}\\
F_j(b_j)  = \left[\begin{matrix}
f_j(b_j)&0 \\
 0&g_j(b_j)
\end{matrix}\right]\qquad 2\le j\le 
d-1
\end{gather*}
and 
$$
F_d(b_d) = \left[\begin{matrix}\alpha^{1/2}f_d(b_d)\\ (1-\alpha)^{1/2} 
  g_d(b_d) \end{matrix}\right] \quad \text{(column  matrix  with 
operator entries).}$$
Then it is easy to check that $\sup\limits_b \|F_j(b)\| < 1$ for
all $1\le j\le 
d$ and hence we obtain 
$$\|\alpha f +
(1-\alpha)g\|_{M_d(I_1,\ldots, I_d)} <1.$$

Moreover, $M_d(I_1,\ldots, I_d)$ is a unital Banach algebra for the
pointwise product, 
i.e.\ for any $f,g$ in $M_d(I_1,\ldots, I_d)$ we have
\begin{equation}\label{eq2.12}
\|f.g\|_{M_d(I_1,\ldots, I_d)} \le \|f\|_{M_d(I_1,\ldots, I_d)} 
\|g\|_{M_d(I_1,\ldots, I_d)}.
\end{equation}
Indeed, if we assume \eqref{sub}   then we have
$$(f.g)(b_1,\ldots, b_d) = (f_1(b_1) \otimes g_1(b_1)) \ldots
(f_d(b_d) \otimes 
g_d(b_d)),$$
and \eqref{eq2.12} follows easily from this. Obviously,
 the
function identically 
 equal to 1 on $I_1\times \ldots\times I_d$ is a unit for
this algebra   and it has   norm 1 in $M_d(I_1,\ldots, I_d)$. 
\begin{exm}\label{exm2.11}
 To illustrate the preceding concepts, we
recover here 
the 
following result from \cite{PyS}:\ Let $G={\bb F}_\infty$ and let ${\cl
W}(1)\subset G$ 
be the subset of all the words of length 1. Then
the indicator function of  ${\cl W}(1)$ is in $ UB(G)$. 
Indeed, we claim that for any bounded function 
$\varphi$ with support in ${\cl W}(1)$ we have
\begin{equation}\label{eq2.14}
\forall~d\ge 1\qquad \|\varphi\|_{M_d(G)} \le
2^d\|\varphi\|_{\ell_\infty(G)}. 
\end{equation}
Thus (by Corollary \ref{cor2.6}) we have $\varphi\in B_c(G)$ for
$c>2$ (actually, this is known for all $c>1$).\hfill\break
 However it can be 
shown that for any such function we have
\begin{equation}\label{eq2.15}
\left(\sum_{t\in G}|\varphi(t)|^2\right)^{1/2} \le 2
\|\varphi\|_{B(G)}. 
\end{equation}
Thus for instance the indicator function of ${\cl W}(1)$ is in
$B_c(G)$ for all 
$c>2$ but not in $B(G)$ (in particular this shows that $G$ is not
unitarizable). 
Note however that since $1_{{\cl W}(1)}$ belongs to $M_d(G)$ for
all $d\ge 1$, 
this does not distinguish the various classes $M_d(G)$ or
$B_c(G)$, but this 
task is completed in \S \ref{sec5}.
\end{exm}
\begin{proof}[Proof of \eqref{eq2.15}]
 First it suffices to prove this
for a finitely supported $\varphi$, with support in ${\cl
W}(1)$. Then
 \eqref{eq2.15} is an immediate consequence of an inequality proved first
by Leinert
\cite{Lei}, and generalized by Haagerup \cite[Lemma 1.4]{H4}:
 Any $\psi$ finitely supported, with support in ${\cl W}(1)$
satisfies
$\|\sum \lambda(t) \psi(t)\|\le 2 (\sum |\psi(t)|^2)^{1/2}$.
The inequality 
\eqref{eq2.15} can be   deduced from this by duality, using the fact
that, if $\varphi$ is finitely supported,  
then $\| \varphi\|_{B(G)}=\sup |\langle  \varphi, \psi \rangle|$
where the
sup runs over all $\psi$ finitely supported on $G$ 
such that $\|\sum \lambda(t) \psi(t)\|\le 1$.
\end{proof}

\begin{proof}[Proof of \eqref{eq2.14}] Let $\varphi$ be a function with
support in ${\cl 
W}(1)$. Consider the set
$$\Omega = \{(t_1,\ldots, t_d)\in G^d\mid t_1t_2\ldots t_d\in
{\cl W}(1)\}.$$
Clearly, when $t_1t_2\ldots t_d$ has length one, it reduces to a
single letter 
(i.e.\ a generator or its inverse). Clearly this letter must
``come'' from 
either $t_1,t_2,\ldots$ or $t_d$. Thus we have   
$$\Omega = \Omega_1 \cup \ldots \cup \Omega_d$$
where $\Omega_j$ is the set of $(t_1,\ldots, t_d)$ in $\Omega$
such that the 
single ``letter'' left after reduction comes from $t_j$.
Hence we have \begin{equation}\label{eq2.151}
 1_{\Omega}=\sum_j  1_{\Omega_j} \prod_{i<j} [1-1_{\Omega_i}]. 
\end{equation} 
For any $\theta$ in $G$, we introduce the operator
$\xi(\theta)\in B(\ell_2(G))$ 
defined as follows:\ Assume $\theta = a_1a_2\ldots a_k$ (reduced
word  where $a_q 
\in {\cl W}(1)$ for all $q$), with $k\ge 1$,  then we set $a_0=a_{k+1}=e$ and
$$\xi(\theta) = \sum^k_{q=1} \varphi(a_q) e_{a_1 \ldots a_{q-1},
(a_{q+1}\ldots 
a_k)^{-1}}$$
where, as usual, $e_{s,t}$ denotes the operator defined by
$e_{s,t}(\delta_t) = 
\delta_s$ and $e_{s,t}(\delta_x) = 0$ whenever $x\ne t$.
Moreover,
if $\theta=e$ (empty word, corresponding to $k=0$), we set
$\xi(\theta) =0$.
 Then it is a simple 
verification that
$$\langle \lambda(t_1) \ldots\lambda(t_{j-1}) \xi(t_j)
\lambda(t_{j+1}) \ldots 
\lambda(t_d) \delta_e,\delta_e\rangle = \varphi(t_1t_2\ldots t_d) 
1_{\{(t_1,\ldots, t_d)\in\Omega_j\}}.$$
A moment of reflection shows that $\|\xi(\theta)\| = 
\sup\limits_q|\varphi(a_q)|$ hence $\sup\limits_{\theta\in
G}\|\xi(\theta)\| \le 
\|\varphi\|_{\ell_\infty(G)}$. This shows with the notation
introduced in Remark 
\ref{rem2.10}, that if we set
$$\Phi_j(t_1,\ldots, t_d) = \varphi(t_1t_2\ldots t_d)
1_{\{(t_1,t_2,\ldots, 
t_d)\in \Omega_j\}}$$
we have
$$\|\Phi_j\|_{M_d(G,\ldots, G)} \le \sup_{t\in G}\|\xi(t)\| \le 
\|\varphi\|_{\ell_\infty(G)},$$
and hence  with $\varphi=1$ identically, we find
$\| 1_{\Omega_j}\|_{M_d(G,\ldots, G)}\le 1$, and
$$\|1- 1_{\Omega_j}\|_{M_d(G,\ldots, G)}\le 2$$
hence by Remark \ref{rem2.10}, \eqref{subadd} and \eqref{eq2.12}, we have
$$\|\varphi\|_{M_d(G)}  = \|\Phi\|_{M_d(G,\ldots, G)}  =
\left\|\sum_j \Phi_j \prod_{i<j} [1-1_{\Omega_i}] \right\|_{M_d(G,\ldots,
G)}
\le  2^d\|\varphi\|_{\ell_\infty(G)}$$ which completes the proof of
\eqref{eq2.14}.
\end{proof}
 
\begin{exm}\label{exm2.12}
 Let $G$ be a free group. 
\begin{itemize}
\item[(i)] Let $\psi_d\colon \ G^d \to \{0,1\}$ be the indicator
function of 
the set formed by all the $d$-tuples $(t_1,\ldots, t_d)$ of
reduced words such 
that $t_i\ne e$ for all $i$ and the product $t_1t_2\ldots t_d$
allows no 
reduction. Then 
$$\|\psi_d\|_{M_d(G,\ldots, G)}\le 5^{d-1}$$
\item[(ii)] A fortiori, for any subsets $I_1\subset G,\ldots,
I_d\subset 
G$, we have
$$\|{\psi_d}_{|I_1\times\cdots\times I_d}\|_{M_d(I_1,\ldots, I_d)}
\le 5^{d-1}.$$
\end{itemize}
\end{exm}

\begin{proof}
Fix $1\le j\le d-1$. Let $A_j$ be the subset of $G^d$ formed of
all 
$(t_1,\ldots, t_d)$ in $G^d$ such that $t_j\ne e$, $t_{j+1}\ne e$
and such that 
 $t_jt_{j+1}$ does reduce, i.e.\ $|t_jt_{j+1}| < |t_j| +
|t_{j+1}|$. Also let 
$B_j = \{t\in G^d\mid t_j\ne e, t_{j+1}\ne e\}$. We will use the
fact that
\begin{equation}\label{eq2.13}
\psi_d = \prod^{d-1}_{j=1} 1_{B_j}-1_{A_j}.
\end{equation}
Observe that for all $t=(t_j)\in G^d$
$$1_{B_j} (t) = (1-1_{\{t_j=e\}})(1-1_{\{t_{j+1}=e\}})$$
and using $1_{\{t_j=e\}} = \langle \delta_{t_j},\delta_e\rangle$,
it is easy to 
deduce from this with \eqref{subadd} and \eqref{eq2.12} that
$$\|1_{B_j}\|_{M_d(G,\ldots, G)}\le 4.$$
Now, for any $x$ in $G$ with $x\ne e$ let us denote by $F(x)$ and
$L(x)$ 
respectively the first and last letter of $x$ (i.e.\ $F(x)$ and
$L(x)$ are 
equal to a generator or the inverse of one). Then it is easy to
check that for 
any $t=(t_j)$ in $G^d$ we have
$$1_{A_j}(t) = \langle\alpha(t_j), \beta(t_{j+1})\rangle$$
where $\alpha(t) = \delta_{L(t)}$ and $\beta(s) =
\delta_{F(s)^{-1}}$ if both 
$|t|>0$ and $|s|>0$ and $\alpha(e) = \beta(e) = 0$. This implies
immediately 
that
$$\|1_{A_j}\|_{M_d(G,\ldots, G)} \le 1,$$
hence $\|1_{B_j}-1_{A_j}\|_{M_d(G,\ldots, G)} \le 5$, and since,
by \eqref{eq2.12},  $M_d(G,\ldots, 
G)$ is a Banach algebra, by \eqref{eq2.13} we obtain 
$$\|\psi_d\|_{M_d(G,\ldots, G)}\le 5^{d-1}.\qquad \qed$$
\renewcommand{\qed}{}\end{proof}

\section{The predual $\pmb{X_d(G)}$ of $\pmb{M_d(G)}$}\label{sec3}

 The definition of the spaces $B_c(G)$ and $M_d(G)$ shows that
they are dual spaces.
There is a natural duality between these spaces
and the group algebra ${\bb C}[G]$ which we view as
the convolution algebra of finitely supported functions on $G$.
Indeed, for any function $f\colon\ G\to {\bb C}$
and any  $g$ in ${\bb C}[G]$, we set
$$<g,f>= \sum_{t\in G} g(t)f(t).$$
Then we define the spaces $X_d(G)$ and $\tilde A_c$
respectively as the completion of ${\bb C}[G]$
for the respective norms
$$\|g\|_{X_d(G)}= \sup \{ |<g,f>| \mid f\in M_d(G), \
\|f\|_{M_d(G)}\le 1
\} $$
and
$$\|g\|_{\tilde A_c}= \sup \{ |<g,f>| \mid f\in M_d(G), \
\|f\|_{B_c(G)}\le 1
\} .$$
Obviously, we can also write
$$\|g\|_{\tilde A_c}= \sup \{ \|\sum g(t) \pi(t)\|\mid
\pi\colon \ G\to B(H),\  |\pi|\le c
\} .$$
This last formula shows that ${\tilde A_c}$ is naturally
equipped with
a Banach algebra structure under convolution: we have
$\|g_1 *g_2\|_{\tilde A_c}\le \|g_1  \|_{\tilde A_c} \| 
 g_2\|_{\tilde A_c}.$

\n
However, the analog for the spaces $X_d(G)$ fails in general.
This was the basic idea
used by Haagerup \cite{H1} to prove that
$M_2({\bb F}_\infty)\not=UB({\bb F}_\infty)$. Indeed, Haagerup used 
spherical functions to show that $X_2({\bb F}_\infty)$ is not
a Banach algebra under convolution (see   Remark \ref{rem3.2} below),
 which implies by the
preceding remarks that $X_2({\bb F}_\infty)\not =\tilde A_c$ for
any
$c$, hence $M_2({\bb F}_\infty)\not =B_c({\bb F}_\infty)$ for
any
$c$, from which $M_2({\bb F}_\infty)\not=UB({\bb F}_\infty)$ follows
easily by Baire's classical theorem.

Note   that, in sharp contrast, for $G={\bb N}$, it is known  that
$X_2(G)$
is a Banach algebra (due to G. Bennett),
 but not an operator algebra (see \cite{P6} for details).

\n Although $X_d(G)$  is not in general
a Banach algebra  under convolution, it satifies the following
property: if $g_1\in X_d(G)$ and $g_2\in X_k(G)$,
then $g_1 *g_2\in X_{d+k}(G)$ and
\begin{equation}\label{eq3.1}
\|g_1 *g_2\|_{ X_{d+k}(G)}\le \|g_1  \|_{X_d(G)} \| 
 g_2\|_{X_k(G)}.
 \end{equation}
Therefore, Haagerup's result in \cite{H1} implies
that $X_2({\bb F}_\infty)\not =X_4({\bb F}_\infty)$
(equivalently $M_2({\bb F}_\infty)\not
=M_4({\bb F}_\infty)$), since otherwise
$X_2({\bb F}_\infty)$ would be a Banach algebra  under convolution.

To verify  \eqref{eq3.1}, we will need an alternate description of
the space  $X_d(G)$, which uses the Haagerup tensor product
and the known results on multilinear cb maps (\cf \cite{ChS,PS}).
These results show that $X_d(G)$  may be identified with 
a quotient (modulo the kernel of the natural product map) of the
Haagerup tensor product
$\ell_1(G)\otimes_h\ldots \otimes_h\ell_1(G)$ of
$d$ copies of $\ell_1(G)$ equipped with its ``maximal operator
space structure".   More explicitly, one can prove
that the space $X_d(G)$ coincides with the space
of all functions $g\colon\ G\to {\bb C}$ for which there is an
element
${\hat g}=\sum_{G^d} {\hat g}(t_1,\ldots ,t_d)  \delta_{t_1}
\otimes\ldots \otimes
\delta_{t_d}$ in 
$\ell_1(G)\otimes_h\ldots \otimes_h\ell_1(G)$ such that
$$\forall t\in G\qquad g(t)=\sum_{t_1\ldots t_d=t} {\hat g}(t_1,\ldots ,t_d)
$$ 
and moreover we have
\begin{equation}\label{eq3.2}
\| g\| _{X_d(G)}= \inf \{\| {\hat g}\|_{
\ell_1(G)\otimes_h\ldots \otimes_h\ell_1(G)}\}.
\end{equation}
In addition the norm of an element ${\hat g}$ in  the space
$\ell_1(G)\otimes_h\ldots \otimes_h\ell_1(G)$ can also be
explicited as follows:  
$$ \| {\hat g} \|_{
\ell_1(G)\otimes_h\ldots \otimes_h\ell_1(G)}= \sup  \{  \left\|
\sum_{G^d} {\hat g}(t_1,\ldots ,t_d) x^1_{t_1}\ldots  x^d_{t_d}
\right\| \} $$
where the supremum runs over all families $(x^1_t)_{t\in
G}$,... , $(x^d_t)_{t\in
G}$ in the unit ball of $B(\ell_2)$. Actually
(by e.g.\ \cite[prop. 6.6]{P8}), the supremum remains the same
 if we restrict it to the case 
when the $d$ families actually coincide with a single
family  $(x_t)_{t\in G}$ in the unit ball of $B(\ell_2)$.

Clearly, 
if ${\hat g}_1\in {\ell_1(G)\otimes_h\ldots \otimes_h\ell_1(G)}$
($d$ times)  and 
${\hat g}_2\in {\ell_1(G)\otimes_h\ldots \otimes_h\ell_1(G)}$
($k$ times)
we have ${\hat g}_1 \otimes {\hat g}_2\in
{\ell_1(G)\otimes_h\ldots \otimes_h\ell_1(G)}$
($d+k$ times) and 
$\|{\hat g}_1 \otimes {\hat g}_2\|\le \|{\hat g}_1  \| \| 
{\hat g}_2\|$. From this,     \eqref{eq3.1}  follows easily using  
\eqref{eq3.2}.

\begin{rk}
 Assume that $ M_d(G)=M_{2d}(G)$.
Then passing to the preduals, $ X_d(G)=X_{2d}(G)$
with equivalent norms. By  \eqref{eq3.1} with $k=d$, this
implies that $ X_d(G)$ is (up to isomorphism)
a Banach algebra under convolution.
Moreover, since  the product
in $ X_d(G)$ is ``induced" by the Haagerup tensor
product,   Blecher's characterization of
operator algebras (see \cite{B} which extends \cite{BRS})
shows that $X_d(G)$ must be (unitally) isomorphic to a (unital)
operator algebra. Combined with Theorem  \ref{thm2.9}, this
implies 
\end{rk}

\begin{thm}\label{thm3.1}
 In the situation of
Theorem  \ref{thm2.9}, the following assertions are
equivalent:
\begin{itemize}
\item[\rm (i)] There is a $\theta\ge 1$ such
that 
$B_{\theta}(G) =
B_c(G)$ for all $c>\theta$.
\item[\rm (ii)]
 There are   $\theta\ge 1$ and  
an integer $d$ such that
$B_{\theta}(G) =M_d(G)$.
\item[\rm (iii)] There is
an integer $d$ such that  $ M_d(G)=M_{2d}(G)$.
\item[\rm (iv)] There is
an integer $d$ such that $ X_d(G)$ is (up to isomorphism)
a unital operator algebra under convolution.
\end{itemize}
\end{thm}

\begin{proof}
 By Theorem  \ref{thm2.9}, (i) implies (ii). By  \eqref{eq2.2} and
 \eqref{eq2.3}, (ii) implies (iii). The
preceding remark shows that  (iii) implies
(iv). Finally,  assume (iv). Then there is a
unital operator algebra
$A\subset B(H)$ and a unital isomorphism $u\colon
\ X_d(G)
\to A$. Let $\theta = \|u\|$ and $K=\|u^{-1}\|$.
Clearly $u$ restricted to the group elements
defines a unital homomorphism
$\pi$ with $|\pi| \le \theta$. By the
very definition of $\|g\|_{\tilde A_\theta}$,
this implies $\|u(g)\|\le \|g\|_{\tilde
A_\theta}$ for all finitely supported $g$,
 hence $\|g\|_{X_d(G)}\le K\|u(g)\|\le
K\|g\|_{\tilde A_\theta}$. Conversely,
we trivially have (see   \eqref{eq2.3}) 
$\|g\|_{\tilde A_\theta}\le \theta^d
\|g\|_{X_d(G)}$. Thus we  obtain  
$X_d(G)={\tilde A_\theta}$, hence by duality
$M_d(G)={B_\theta(G)}$, which
(recalling the basic inclusions  \eqref{eq2.2} and  \eqref{eq2.3}))
implies (i). 
\end{proof}

\begin{rem}\label{rem3.2}
 Haagerup's proof in \cite{H1} that
$M_2({\bb F}_\infty)\ne 
UB({\bb F}_\infty)$ can be outlined as follows. Let $G = {\bb F}_n$ with
$2\le
n<\infty$.  Assume $M_2(G) = UB(G)$. 
\end{rem}

\n Step 1:\ By Baire's theorem, there exists $c>1$ such that
$M_2(G) = B_c(G)$ 
with equivalent norms.

\n Step 2:\ This implies that $X_2(G)$ is a Banach algebra under
convolution 
(because the predual of $B_c(G)$ is clearly an operator algebra,
see Theorem~\ref{thm3.1} 
above). Hence, there is $C>0$ such that for all $f,g$ finitely
supported   we have
$\|f*g\|_{X_2(G)}\le C\|f\|_{X_2(G)} \|g\|_{X_2(G)}$.

\n Step 3:\ By an averaging argument, the radial projection $f\to
f_R$ defined 
by
$$f_R(t) = \sum_{s\colon \ |s|=|t|} f(s) \cdot
[\text{card}\{s\mid |s| = 
|t|\}]^{-1}$$
is bounded on $M_2(G)$ so that $\|f_R\|_{M_2(G)}\le
\|f\|_{M_2(G)}$ for any $f$ 
in $M_2(G)$.

\n Step 4:\ Let $\varphi_z$ be the spherical function on $G$
equal to $z$ on 
words of length 1 (cf.\ e.g.\ \cite{FTP}). This means that
$\varphi_z(t) = \varphi(|t|)$ 
where $\varphi$ is determined   inductively by:\ $\varphi(0) =
1$, $\varphi(1) = 
z$ and
$$\varphi(k+1) = \frac{2n}{2n-1} \varphi(1) \varphi(k) - \frac{1}{2n-1} 
\varphi(k-1)$$
for all $k\ge 2$.
The spherical property of $\varphi_z$ implies that for any
finitely supported 
radial function $f$ we have 
$\varphi_z*f  = \langle\varphi_z,f\rangle\ \varphi_z,$ and hence
if $g$ is another
finitely supported 
radial function,   we have
$$ \langle\varphi_z,f*g\rangle=\langle\varphi_z,f\rangle
\langle\varphi_z,g\rangle.$$
Moreover, if $|z|<1$ then $\varphi_z\in M_2(G)$ (actually
$\varphi_z\in UB(G)$, see \cite{MPSZ}).
Thus, in short, although Step 5 below says it is unbounded,
  $\|\varphi_z\|_{M_2(G)}$ is finite whenever $|z|<1$.
Therefore, if $|z|<1$, $f\to \langle\varphi_z,f\rangle$ defines a
continuous multiplicative unital functional
on the Banach subalgebra which is the closure of the set
of  finitely supported radial functions in $X_2(G)$.
Clearly, this implies that $\langle\varphi_z,f\rangle$ is in the
spectrum of
$f$, hence its modulus is  majorized by
the spectral radius of $f$ in the latter Banach algebra, 
and this is $\le C\|f\|_{X_2(G)}$ by Step 2.
 Thus we obtain for $f$ radial    $ |   \langle\varphi_z,f\rangle
|  \le C\|f\|_{X_2(G)}$.
\n Now, for $f$ finitely supported but not necessarily radial, we
have
\begin{align*}
|\langle\varphi_z,f\rangle| &=
|\langle\varphi_z,f_R\rangle|\\\
&\le C\|f_R\|_{X_2(G)}\\
\intertext{hence by Step 3}
 &\le C\|f\|_{X_2(G)}.
 \end{align*}
This implies $\|\varphi_z\|_{M_2(G)}\le C$. But this contradicts
the next and final step
proved in \cite{H1}:
\medskip

\n Step 5:\ $\sup\limits_{|z|<1} \|\varphi_z\|_{M_2(G)} =\infty$.

\n  It is not 
clear to us how to extend this argument to $M_d$ in place of
$M_2$. The analogue 
of Step 3 is not clear to us (but seems likely to be true).
 Moreover, note that the analogue of Step 2 would 
require assuming $X_d(G) = X_{2d}(G)$, so it would seem that the 
argument would lead, at best, to $M_d(G) \ne M_{2d}(G)$.

\section{The $B(H)$-valued case}\label{sec4}

 Up to now, we have mainly concentrated on
properties  of spaces of coefficients or of
analogous spaces of complex valued functions on
$G$. We now turn to the more general $B(H)$-valued
case which is entirely similar to the preceding
treatment (corresponding to $\dim(H)=1$).
 More generally, for any u.b.\ representation
$\pi\colon \ G\to B(H)$ let us  define
$$\text{Sim}(\pi) = \inf\{\|{S}\|\ \|{S}^{-1}\| \mid
{S}^{-1}\pi(\cdot){S} \text{ is  a unitary
representation}\},$$ 
and let
$$\alpha(G) = \inf\{\alpha\ge 0 \mid \exists K \ \forall
\pi\colon \ G\to B(H) 
\text{ u.b. } \text{Sim}(\pi) \le K|\pi|^\alpha\}.$$
Then again the same phenomenon arises:

\begin{thm}\label{thm4.1}
 Let $G$ be a discrete group.
If $G$ is unitarizable then 
$\alpha(G) < \infty$. Moreover $\alpha(G) \in {\bb N}$.
\end{thm}

We will now explain what replaces $d_1(G)$ in this case.

First, we need to generalize the space $B(G)$, from complex
values to operator 
values. Let $H$ be a Hilbert space and let $G$ be a semi-group
with unit. We denote by $B_c(G;H)$  
the space of all functions $f\colon \ G\to B(H)$ for which there
are a 
 u.b.\  unital homomorphism $\pi\colon \ G\to B(H_\pi)$ with 
$|\pi|\le c$ and operators $\xi\colon \ H_\pi\to H$ and
$\eta\colon \ H\to 
H_\pi$ such that
$$f(t) = \xi\pi(t)\eta.\leqno \forall~t\in G$$
We define
$$\|f\|_{B_c(G;H)} = \inf\{\|\xi\|\ \|\eta\|\}$$
  where the infimum runs 
over all possible such representations. 

\n Here again, in the group case
we will denote $B_1(G;H)$  simply by $B(G;H)$ to emphasize
that  $|\pi|\le 1$ means that $\pi$ is a unitary representation.

Similarly, we denote by $M_d(G;H)$ the space of functions
$f\colon \ G\to B(H)$ 
for which there are bounded functions $\xi_i\colon \ G\to B(H_i,
H_{i-1})$, 
$1\le   i\le d$, with $H_0=H_d=H$, $H_i$ Hilbert such that
$$\forall(t_1,\ldots, t_d)\in G^d\qquad f(t_1t_2\ldots t_d) = 
\xi_1(t_1)\xi_2(t_2)\ldots
\xi_d(t_d).
$$
We equip this space with the norm
$$\|f\|_{M_d(G;H)} = \inf\{\sup_{t_1\in G} \|\xi_1(t_1)\|\ldots
\sup_{t_d\in G} 
\|\xi_d(t_d)\|\}.$$
Note that there is also an obvious
$B(H)$-valued generalization of the spaces
$M_d(I_1,\ldots ,I_d)$ introduced in Remark \ref{rem2.10} above.
Let us denote it by $M_d(I_1,\ldots ,I_d;H)$.
Then,  as before,  for any $\varphi\colon\ G \to B(H)$,
let $\Phi\colon\ G^d \to B(H)$ be defined by
$\Phi(t_1,\ldots ,t_d)=\phi(t_1t_2\ldots t_d)$. Then we have
$$\|\varphi\|_{M_d(G;H)} = \|\Phi\|_{M_d(G,\ldots, G;H)}.$$

Clearly, when $\dim(H) =1$, we recover the previous spaces
$B(G)$, $B_c(G)$ and 
$M_d(G)$. The following extensions of the previous results can be
proved:

\begin{thm}\label{thm4.2}
 Consider a function
$f\colon \ G\to B(H)$. Then $f\in  B_1(G;H)$ iff
$f\in \bigcap\limits_{m\ge 1} M_m(G;H)$ and
$\sup\limits_m 
\|f\|_{M_m(G;H)} <\infty$.
Moreover we have 
$$\|f\|_{B_1(G;H)} =
\sup\limits_{m\ge 1} 
\|f\|_{M_m(G;H)}.$$
 On the other hand, $f\in
\bigcup\limits_{c>1} B_c(G;H)$ iff  
$f\in \bigcap\limits_m M_m(G;H)$ and $\sup\limits_{m\ge 1} 
\|f\|^{1/m}_{M_m(G;H)} < \infty$.
In addition 
$$\limsup_{m\to \infty} \|f\|^{1/m}_{M_m(G;H)}=\inf\{c\ | \ f\in
B_c(G;H)\}.$$
\end{thm}  

\begin{note}
 Let us denote by $d(G)$ the smallest $d$ such
that
$$B(G;\ell_2) = M_d(G;\ell_2).$$
Then we have:
\end{note}

\begin{thm}\label{thm4.3}
 For any unitarizable group,
we have
$$\alpha(G) = d(G).$$
\end{thm}

\n {\bf Warning.} Unfortunately no example is
known of $G$ with $3\le \alpha(G)< 
\infty$.

\begin{rk}
 Theorem \ref{thm4.3} (with Theorems \ref{thm2.3} and \ref{thm0.1}) shows that
$\alpha(G)<3$ iff $G$ is amenable.
\end{rk}

We now turn to the $B(H)$-valued 
variant of the space
$X_d(G)$.  Here we will   
use explicitly the Haagerup tensor product for
operator spaces. We refer the
reader to \cite{ChS,PS} for more on this notion.

\n Let $H$ be an infinite dimensional  Hilbert space.
We
denote by
${K(H)}$ the space of compact operators on $H$.

\n Let $E_1,E_2$ be operator spaces.
Let $x_1\in
{K(H)}\otimes E_1$, $x_2\in
{K(H)}\otimes E_2$. We will denote by
$(x_1,x_2) \to x_1\odot x_2$ the bilinear
mapping from
$({K(H)}\otimes E_1) \times ({K(H)}\otimes E_2)$ to
${{K(H)}}\otimes (E_1\otimes E_2)$ which
is defined on rank one tensors by
$$(k_1\otimes e_1) \odot (k_2 \otimes e_2) = (k_1k_2) \otimes
(e_1\otimes e_2).$$

The Haagerup tensor product $E_1\otimes_h E_2$
 can be characterized as the unique operator space
which  is a completion of
the algebraic tensor product and is such that
for any $x\in K(H)\otimes [E_1\otimes_h E_2]$
we have
$$\|x\|_{\min} = \inf \{
\|x_1\|_{\min}\|x_2\|_{\min}\}$$
where the infimum runs over all 
factorization of the form
$$x=x_1\odot x_2$$
with $x_1\in K(H)\otimes  E_1$ 
and $x_2 \in
K(H)\otimes  E_2$.

By definition of the Haagerup tensor product,
$(x_1,x_2) \to x_1\odot x_2$   extends to a
contractive bilinear  mapping from
$({K(H)}\otimes_{\min} E_1) \times
({K(H)}\otimes_{\min} E_2)$ to
${{K(H)}}\otimes_{\min} [E_1\otimes_h E_2]$. We will
still denote by
$(x_1,x_2) \to x_1\odot x_2$
this extension, and similarly for  $d$-fold tensor products.

The space $X_d(G;H)$ is defined as  the space
of all functions $g\colon\ G\to K(H)$ for which there is an
element
${\hat g}=\sum_{G^d} {\hat g}(t_1,\ldots ,t_d) \otimes
\delta_{t_1}
\otimes\ldots \otimes
\delta_{t_d}$ in 
$K(H)\otimes_{\min}[\ell_1(G)\otimes_h\ldots \otimes_h\ell_1(G)]$
such that
\[
\forall t\in G\qquad  g(t)=\sum_{t_1\ldots t_d=t} {\hat g}(t_1,\ldots ,t_d)
\]
 and moreover we have
\begin{equation}\label{eq4.1}
\| g\| _{X_d(G;H)}= \inf \{\| {\hat g}\|_{K(H)\otimes_{\min}
[\ell_1(G)\otimes_h\ldots \otimes_h\ell_1(G)]}\}.
\end{equation}

In addition the norm of an element ${\hat g}$ in  the space
$K(H)\otimes_{\min}[\ell_1(G)\otimes_h\ldots \otimes_h\ell_1(G)]$
can also be explicited as follows:  
$$ \| {\hat g}\|_{K(H)\otimes_{\min}[
\ell_1(G)\otimes_h\ldots \otimes_h\ell_1(G)]}= \sup \left \{ \|
\sum_{G^d} {\hat g}(t_1,\ldots ,t_d)\otimes  x^1_{t_1}\ldots 
x^d_{t_d}
\|_{B(H\otimes \ell_2)}
\right\}$$ 
where the supremum runs over all families $(x^1_t)_{t\in
G}$,... , $(x^d_t)_{t\in
G}$ in the unit ball of $B(\ell_2)$. Here again
(by e.g.\ \cite[prop. 6.6]{P8}) the supremum is
the same if we restrict the supremum to the case 
when the $d$ families are all equal to a single  
family  $(x_t)_{t\in G}$ in the unit ball of $B(\ell_2)$.
 
By definition of the Haagerup tensor product,
we have also
\begin{equation}\label{eq4.2}
 \| {\hat g}\|_{K(H)\otimes_{\min}[
\ell_1(G)\otimes_h\ldots \otimes_h\ell_1(G)]}= 
\inf \left \{  \|g_1\|_{K(H)\otimes_{\min}
\ell_1(G)}  \cdots
\|g_d\|_{K(H)\otimes_{\min}
\ell_1(G)} \right\}
\end{equation}
where the infimum runs over all
factorizations of ${\hat g}$ of the form
\begin{equation}\label{eq4.3}
{\hat g}= g_1\odot g_2\odot \cdots \odot g_d,
\end{equation}
with $g_1,g_2,\cdots,g_d \in {K(H)\otimes_{\min}
\ell_1(G)}$. Equivalently,  \eqref{eq4.3} means that  for all
$(t_i)$ in
$G^d$ we have
$$ {\hat g}(t_1,\ldots ,t_d)=g_1(t_1) g_2(t_2)\cdots g_d(t_d),
 $$
where the product is in $K(H)$.

From this,    we deduce that 
\begin{equation}\label{eq4.4}
\| g\| _{X_d(G;H)}= \inf \left \{  \|g_1\|_{K(H)\otimes_{\min}
\ell_1(G)}  \cdots
\|g_d\|_{K(H)\otimes_{\min}
\ell_1(G)} \right\},
\end{equation}
where the infimum runs over all factorizations  of $g$  
(as a generalized $d$-fold
convolution) of the form
$$g(t)= \sum_{t_1\ldots t_d=t}g_1(t_1) g_2(t_2)\cdots
g_d(t_d)\qquad (g_i\in {K(H)\otimes_{\min}
\ell_1(G)}).$$

In particular, \eqref{eq4.4}   implies  that
for any integers $d,k$, we have
\begin{equation}\label{eq4.5}
\| g\| _{X_{d+k}(G;H)}=\inf \left \{ \|x\|_{\min}\|y\|_{\min}
\right\},
\end{equation}
where the infimum runs
over all
pairs
$x\in {X_{d}(G;H)}$  $y\in {X_{k}(G;H)}$ such that\hfill\break
\centerline{$g(t)= \sum_{t_1 t_2=t} x(t_1) y(t_2)$  
{($K(H)$-valued convolution)}.}

 The next statement is the $B(H)$-valued
analogue of Theorem  \ref{thm3.1}.

\begin{thm}[\cite{P4}]\label{thm4.4}
Let $G$ be a
semigroup with unit. Let $H=\ell_2$.
 The following assertions are
equivalent:
\begin{itemize}
\item[\rm (i)] There is a $\theta\ge 1$ such
that 
$B_{\theta}(G;H) =
B_c(G;H)$ 
for all $c>\theta$.
\item[\rm (i)$'$] There is a $\theta\ge 1$ such
that 
$B_{\theta}(G;H) =
B_c(G;H)$ 
for some $c>\theta$.
\item[\rm (ii)]
 There are   $\theta\ge 1$ and  
an integer $d$ such that
$B_{\theta}(G;H) =M_d(G;H)$.
\item[\rm (iii)] There is
an integer $d$ such that  $ M_d(G;H)=M_{d+1}(G;H)$.
\item[\rm (iv)] There is
an integer $d$ such that $ X_d(G)$ is (up to
complete isomorphism) a unital operator algebra under
convolution.
\end{itemize}
\end{thm}

\begin{proof}
In \cite{P3} the above Key Lemma \ref{lem2.4} is actually proved in the
$B(H)$-valued case. Therefore all the preceding statements
numbered between 2.5 and 2.9 remain valid in 
the
$B(H)$-valued case. Thus
exactly the same rasoning as  for Theorem  \ref{thm3.1} yields
the equivalence of (i), (i)' and (ii). Clearly (ii) implies
(iii).

\n   Assume (iii).
 Then, passing to the preduals (here we mean the
preduals of $M_d(G)$ and $M_{d+1}(G)$ 
in the operator space sense),
we find  $ X_d(G;H)=X_{d+1}(G;H)$. This
implies  $ X_{d+1}(G;H)=X_{d+2}(G;H)$. Indeed,
by
\eqref{eq4.5} 
  for any $g$ in
the open unit ball of $ X_{d+1}(G;H)$, we can find
$x$ in the unit ball of $ X_{d}(G;H)$
and $y$ in the unit ball of $ X_{1}(G;H)$ such that
$$g(t)=\sum_{t_1t_2=t} x(t_1) y(t_2) .$$
Now since  $ X_d(G;H)=X_{d+1}(G;H)$,  $x\in  X_{d+1}(G;H)$
hence  by    \eqref{eq4.5} $g$ must be in
$X_{d+2}(G;H)$. Now from
$ X_{d+1}(G;H)=X_{d+2}(G;H)$ we deduce
$ X_{d+2}(G;H)=X_{d+3}(G;H)$, and so on..,
so that we must have $ X_d(G;H)=X_{2d}(G;H)$, or equivalently
$ X_d(G)=X_{2d}(G)$   completely isomorphically.
Note that, by   \eqref{eq4.4} or  \eqref{eq4.5}, the
convolution product defines a completely contractive
linear map  $p$ from 
$   X_d(G)\otimes_h  X_d(G)$  to $X_{2d}(G)$, hence
 since $ X_d(G)=X_{2d}(G)$   completely isomorphically,
$p$ is c.b.\  from $X_d(G)\otimes_h  X_d(G)$  to $X_{d}(G)$,
which implies by Blecher's result in \cite{B}
that $X_{d}(G)$ is completely isomorphic to an operator algebra.
This proves that (iii) implies (iv).

\n Finally, assume (iv). Then, there are a unital
subalgebra
$A\subset B({\cl H})$ and   
a unital homorphism $u\colon\ X_d(G) \to A$ which is
also a complete isomorphism. Let $\theta=
\|u\|_{cb}$ and
$C= \|u^{-1}\|_{cb}$. Let $\pi(t)=u(\delta_t)$ ($t\in G$). 
Then $\pi$ is a u.b.\ representation of $G$ with $|\pi| \le
\theta$. By the maximality of $\tilde A_\theta$, 
for any $x\in {\bb C}[G]$, we must have
$$ \|u(x)\| \le \|x\|_{\tilde A_\theta} ,$$
hence 
$ \|x\|_{X_d(G)} \le C\|x\|_{\tilde A_\theta}.$
By duality, this implies
that for all $\varphi$ in $M_d(G)$ we have 
$$\|\varphi\|_{B_{\theta}(G)} \le C \|\varphi\|_{M_d(G)}.$$
Moreover, the same arguments with coefficients in
$B(H)$ yield the c.b.\ version of
this, so that we obtain, for all $\varphi$ in $M_d(G;H)$
$$\|\varphi\|_{B_{\theta}(G;H)} \le C \|\varphi\|_{M_d(G;H)}.$$
Thus we obtain (ii) and hence also (i), establishing 
(iv) $\Rightarrow$ (i).
\end{proof}

\begin{rk}
 The preceding argument shows that
(iii) and (iv) are equivalent for the same $d$.
\end{rk}

\section{A case study: The free groups}\label{sec5}

We wish to prove here the following.

\begin{thm}\label{thm5.1}
 For any $d\ge 2$,
$$M_d({\bb F}_\infty) \ne M_{d-1}({\bb F}_\infty).$$
More precisely let $\{g_1,g_2,\ldots\}$ be the
free generators of ${\bf F}_\infty$,  and for any $n$
let $W_{d,n}$ be the subset of ${\bb F}_\infty$ formed
of all the  words $w$ (of length $d$) of the form
$w = g_{i_1}g_{i_2}\ldots g_{i_d}$ with 
$1\le i_j\le n$ for any $1\le j\le d$. Then, for any $n$, there
is a function 
$f_{d,n}\colon \ {\bb F}_\infty\to {\bb C}$ supported on $W_{d,n}$ and 
unimodular
on 
$W_{d,n}$ such that
$$n^{\frac{d-1}{2}} \le \|f_{d,n}\|_{M_d({\bb F}_\infty)}\quad \text{and}\quad 
\|f\|_{M_{d-1}({\bb F}_\infty)}\le C(d)n^{\frac{d-2}{2}},$$
where $C(d)$ is a constant depending only on $d$.
\end{thm}

Let 
$$ {UB}(G)=\bigcup_{c>1}
B_c(G).$$
Since we have obviously inclusions ${UB}(G) \subset M_d(G)
\subset M_{d-1}(G)$ 
for any group $G$, this implies

\begin{cor}\label{cor5.2}
 For any $d\ge 1$,
$$M_d({\bb F}_\infty) \ne {UB}({\bb F}_\infty).$$
\end{cor}

For $d=2$ this is the main result of \cite{H1}. Note however that 
Theorem \ref{thm5.1} yields a function
$f$ supported in the words of length $3$ that is in $M_{2}(G)$
but not in $M_{3}(G)$ 
and hence not in ${UB}(G)$. It is easy to see that $3$ is 
minimal here, i.e.\ any function
supported in the words of length $2$ that is in $M_{2}(G)$ must
be 
in ${UB}(G)$ (see Proposition \ref{pro5.8} below).

  Let $I_1,\ldots, I_d$ be arbitrary sets. 
Recall the notation from Remark \ref{rem2.10}:\  We   denote 
by 
$M_d(I_1,\ldots, I_d)$ the space of all functions $f\colon \ 
I_1\times\cdots\times I_d\to {\bb C}$ for which there are bounded
functions 
$f_i$
$$f_i\colon \ I_i \to B(H_i,H_{i-1})\quad
(\text{here}\  H_i \ \text{are Hilbert spaces} \text{with } H_d =
H_0 = {\bb C}) \text{  such
that}$$
$$f(b_1,\ldots, b_d) = f_1(b_1)\ldots f_d(b_d)\qquad
\forall~b_i\in I_i.$$
We equip this space with the norm
$$\|f\| = \inf\left\{\prod^d_{i=1}
 \sup_{b\in I_i} \|f_i(b)\|\right\}$$
where the infimum runs over all possible such factorizations.

Let $J_i \subset I_i$ be arbitrary subsets. Note that we
obviously have
$$\|f_{|J_1\times\cdots\times J_d}\|_{M_d(J_1,\ldots, J_d)} \le
\|f\|_{M_d(I_1,\ldots, I_d)}.$$
Moreover, for any function $g\colon \ J_1\times\cdots\times
J_d\to {\bb C}$ let 
$\tilde g\colon \ I_1\times\cdots\times I_d\to {\bb C}$ be the
extension of 
$\tilde g$ equal to zero outside $J_1\times\cdots\times J_d$.
Then it is easy to 
check that
\begin{equation}\label{eq5.1}
\|\tilde g\|_{M_d(I_1,\ldots, I_d)} = \|g\|_{M_d(J_1,\ldots,
J_d)}.
\end{equation}

We will relate these spaces to $M_d({\bb F}_\infty)$ via the following
observation. 
Given a function $\varphi\colon \ {\bb F}_\infty\to {\bb C}$ supported by
$W_{d,n}$, we 
can define $f\colon \ [1,\ldots, n]^d \to {\bb C}$ by 
$$\forall i_j\in [1,\ldots, n]\qquad\qquad f(i_1,i_2,\ldots, i_d)
= 
\varphi(g_{i_1}g_{i_2}\ldots g_{i_d}).$$
We have then obviously if $I = [1,\ldots, n]$
\begin{equation}\label{eq5.2}
\|f\|_{M_d(I,\ldots, I)} \le
\|\varphi\|_{M_d({\bb F}_\infty)}.
\end{equation}
The main idea for the proof of Theorem \ref{thm5.1} is to compare 
$\|\varphi\|_{M_{d-1}({\bb F}_\infty)}$ with certain norms of $f$ of the form 
$M_{d-1}(I_1,\ldots,I_{d-1})$ when $f$ is viewed as depending on
less than $d$ 
variables, by blocking together certain variables, so that 
$I_1 = I^{p_1}$, $I_2 = I^{p_2},\ldots$ with $p_1+p_2 +\cdots+
p_{d-1} = d$.

\begin{rk}
 With the notation used in operator space theory,
the space 
$M_d(I_1,\ldots, I_d)$ can be identified with the dual of the
Haagerup tensor 
product $\ell_1(I_1) \otimes_h\cdots\otimes _h \ell_1(I_d)$,
where the 
spaces $ \ell_1(I_j) $
are equipped (as usual) with their maximal operator space
structure
in the sense of e.g.\ \cite{ER2} or \cite{P8}.
\end{rk}

Consider now a partition 
$\pi = (\alpha_1,\ldots, \alpha_k)$ of
$[1,\ldots,  d]$ into disjoint intervals
(=``blocks") with
$k<d$, so that at least one
$\alpha_i$ has $|\alpha_i|>1$. Let 
$I(\alpha_i) = \prod\limits_{q\in \alpha_i} I_q$.

We have a natural mapping from $M_d(I_1,\ldots, I_d)$ to 
$M_k(I(\alpha_1),\ldots, I(\alpha_k))$ associated to the
canonical 
identification
$$I_1 \times\cdots \times I_d = I(\alpha_1)\times\cdots\times 
I(\alpha_k).$$
It is easy to check that this mapping is contractive. For
simplicity of 
notation, we denote
$$M(\pi) = M_k(I(\alpha_1),\ldots, I(\alpha_k)).$$
Moreover, it is useful to observe that if $\pi'$ 
is another partition of 
$[1,\ldots, d]$ that is finer than $\pi$ 
(i.e.\ such that every block of $\pi$ 
is a union of certain blocks of $\pi'$),
 then we have $M(\pi') \subset M(\pi)$ 
and for any $f\colon \ [1,\ldots, n]^d \to {\bb C}$
\begin{equation}\label{eq5.3}
\|f\|_{M(\pi)} \le \|f\|_{M(\pi')}.
\end{equation}
Note however that since the set of all partitions
 is only {\sl partially\/} ordered (and 
not totally ordered), the intersection 
$\bigcap\limits_\pi M(\pi)$ over all 
partitions with $k$ blocks does not reduce 
to one of the $M(\pi)$. We equip this 
intersection ${\bigcap\limits_\pi M(\pi)}$ with
 its natural norm, namely :
$$\|f\| = \max_\pi \|f\|_{M(\pi)},$$
where the
 maximum runs over all $\pi$ with  at most $d-1$ blocks.

The main point in the proof of Theorem \ref{thm5.1}
is the following.

\begin{lem}\label{lem5.3}
 Assuming $I_1,\ldots, I_d$ are infinite
sets, then for any 
$d>1$ the natural mapping
$$\Phi\colon \ M_d(I_1,\ldots, I_d) \longrightarrow \bigcap_\pi
M(\pi)$$
is not an isomorphism.
\end{lem}

To prove this, we will use two more lemmas.

\begin{lem}\label{lem5.4}
 For any $f$ in $M_d(I_1,\ldots, I_d)$ and
any fixed  $j$ in 
$[1,2,\ldots, d]$ we have
$$\|f\|_{M_d(I_1,\ldots, I_d)} \le \sup_{I_1\times\cdots\times
I_d} |f| \cdot 
\left[\prod_{m\ne j} |I_m|\right]^{1/2}.$$
\end{lem}

\begin{proof}
 Let us write
$$f(i_1,\ldots, i_d)= f(a,i_j,b).$$
Then for each fixed $j$, the matrix $(f(a,i_j,b))_{a,b}$ defines
an operator 
$\xi(i_j)$ from $H =  \ell_2(I_{j+1}) \otimes \cdots \otimes
\ell_2(I_d)$ to $K 
= \ell_2(I_1) \otimes\cdots\otimes \ell_2(I_{j-1})$ and its norm
can be 
majorized as follows (observe that an $n\times m$ matrix
$(a_{ij})$ has norm 
bounded by $\sup\limits_{ij} |a_{ij}| \sqrt n~\sqrt m$)
$$\|\xi(i_j)\|\le \sup|f| \cdot \left(\prod_{m\ne
j} |I_m|\right)^{1/2}.$$ 
Moreover we have
$$f(i_1,\ldots, i_d) = a_1(i_1) \ldots a_j(i_{j-1}) \xi(i_j)
b_{j+1}(i_{j+1}) 
\ldots b_d(i_d)$$
with $a_m\colon \ K\to K$ and $b_m \colon \ H\to H$ defined by
$a_m(i_m) =  
1\otimes\cdots 1\otimes e_{1i_m} \otimes 1 \cdots \otimes 1$ and
$b_m(i_m) = 1 
\otimes \cdots 1 \otimes e_{i_m1} \otimes 1 \cdots\otimes 1$,
where the middle 
term is in the place of index $m$. 
\end{proof}

\begin{lem}[Marius Junge]\label{lem5.5}
 Assume
$|I_1| = 
\cdots = |I_d| = n$. Then the natural identity map
 from $\ell_\infty(n^d)$
 to $M_d(I_1,\ldots, I_d)$
has norm $n^{\frac{d-1}{2}}.$ Equivalently,
we have
\begin{equation}\label{eq5.4}
\sup \left\{ \left\| \sum z_{i_1 i_2\ldots i_d}\ 
e_{i_1}\otimes \cdots  \otimes e_{i_d}
\right\|_{M_d(I_1,\ldots, I_d)} \ \mid\ 
|z_{i_1 i_2\ldots i_d}|\le 1 \right\} =n^{\frac{d-1}{2}}.
\end{equation}
\end{lem}

\begin{proof}
 I am grateful to Marius Junge for kindly
providing this lemma in answer to a question of
mine for $d=3$.   Let $C$ be the left side of
\eqref{eq5.4}.   The fact that $C\le n^{\frac{d-1}{2}}$ follows from Lemma 
\ref{lem5.4}.
The main point is the converse. 
To prove this, 
consider the function
$$\psi(i_1,\ldots, i_d) = a_{i_1i_2} a_{i_2i_3}\ldots
a_{i_{d-1}i_d}$$
where $(a_{ij})$ is an $n\times n$ unitary matrix with $|a_{ij}|
= n^{-1/2}$. 
Let $\xi_{i}(x)=\langle x,e_i  \rangle $.
Then we can write
$$\psi(i_1,\ldots, i_d) = \langle a e_{i_2i_2} ae _{i_3i_3}
\ldots e_{i_{d-1}i_{d-1}} ae_{i_d}, 
e_{i_1}\rangle  =
\xi_{i_1}(a e_{i_2i_2}ae _{i_3i_3} \ldots e_{i_{d-1}i_{d-1}}
ae_{i_d})
$$
where $a$ appears $d-1$ times, from which it follows   that
\begin{equation}\label{eq5.5}
\|\psi\|_{M_d(I_1,\ldots, I_d)^*} \le n.
\end{equation}
Indeed, if $\|f\|_{M_d(I_1,\ldots, I_d)}<1$. We may assume
$f(b_1,\ldots, b_d) = f_1(b_1)\ldots f_d(b_d)\qquad
\forall~b_i\in I_i$
with $\|f_i(b_i)\|<1$ ($b_i\in I_i$) for all $i$.
Then we have
$$\sum \psi(i_1,\ldots, i_d) f(i_1,\ldots, i_d)= 
[\sum_{i_1} \xi_{i_1}\otimes f_1(i_1)] [a\otimes I] [\sum_{i_2}
e_{i_2i_2}\otimes f_2(i_2)]\ldots $$
hence
$$|\sum \psi(i_1,\ldots, i_d) f(i_1,\ldots, i_d)|\le  \|\sum
\xi_{i_1}\otimes f_1(i_1)\| \|a\|^{d-1}
 \|\sum  f_d(i_d)  \otimes e_{i_d}\|$$
$$\le \sqrt{n} \|a\|^{d-1} \sqrt{n}\le n,$$
which establishes \eqref{eq5.5}
Therefore we must have
$$\sum|\psi(i_1,\ldots, i_d)| \le Cn,$$
whence
$$n^d(n^{-1/2})^{d-1} \le Cn$$
which yields $C\ge n^{\frac{d-1}{2}}$ as announced.

\n Note:\ The preceding proof uses implicitly
ideas from operator space theory namely 
the identity 
$M_d(I_1,\ldots, I_d)=\ell_\infty^n
\otimes_h\cdots \otimes_h \ell_\infty^n$ ($d$
times), for which we refer to e.g. \cite{ER2} or \cite{P8}.
\end{proof}

\begin{rk}
 Let
$$\vp(p,q) = \exp\{ipq/n\},$$
and let $a_{pq} = \vp(p,q)n^{-1/2}$. Thus, the
$n\times n$ unitary matrix $a =  (a_{pq})$
represents the Fourier transform on the group
${\bb Z}/n{\bb Z}$. Let
$$F_{d,n}(i_1,\ldots, i_d) = \vp(i_1,i_2) \vp(i_2,i_3)\ldots
\vp(i_{d-1},i_d).$$
Then the preceding proof yields
$$\|F_{d,n}\|_{M_d(I_1,\ldots, I_d)} = n^{\frac{d-1}{2}}.$$
\end{rk}

\begin{proof}[Proof of Lemma \ref{lem5.3}]
 By \eqref{eq5.1},
 it suffices to show that if $|I_1| = 
\cdots = |I_d| = n$ then $\|\Phi^{-1}\| \ge \sqrt
n$  for all $n$.
Thus we now assume
$|I_1| =\cdots= |I_d| = n$ throughout this proof.
By  Lemma~\ref{lem5.4}, for any $\pi$ we have
$\|\Phi^{-1}\colon \ \ell_\infty(n^d) \to 
M(\pi)\|\le n^{\frac{d-2}{2}}$. (Indeed, we can
choose $\alpha_j$ with $|\alpha_j| 
\ge 2$, hence $|I(\alpha_j)| \ge n^2$.)
It follows that
$$\left\|\Phi^{-1}\colon \ \ell_\infty(n^d) 
\longrightarrow \bigcap_\pi 
M(\pi)\right\| \le n^{\frac{d-2}{2}}.$$
Note that for the mapping underlying $\Phi^{-1}$
we have 
$$\|  \ell_\infty(n^d) 
\to M_d(I_1,\ldots, I_d)\|\le 
\|   \ell_\infty(n^d) 
\longrightarrow \bigcap_\pi 
M(\pi)\|\times 
\|    \bigcap_\pi 
M(\pi) 
\longrightarrow  
M_d(I_1,\ldots, I_d)\|$$
Thus  the above estimate together with Lemma \ref{lem5.5}
implies
$$\sqrt n\le \|\Phi^{-1}\colon \ \bigcap_\pi 
M(\pi) 
\longrightarrow  
M_d(I_1,\ldots, I_d)\|.\qquad\qed$$
\renewcommand{\qed}{}\end{proof}

Let $G = {\bb F}_n$ with $2\le n\le \infty$ 
and let $I$ denote the set of generators 
of $G$. Let $W_d$ be 
the set of all elements of $G$ which 
are a product of exactly $d$ generators. 
Let $F\colon \ I^d\to {\bb C}$ be a function 
and let $f\colon \ G\to {\bb C}$ be 
the function defined on $W_d$ by
\begin{equation}\label{eq5.6}
f(i_1i_2\ldots i_d) = F(i_1,i_2,\ldots, i_d)
\end{equation}
and equal to zero outside $W_d$.

\begin{lem}\label{lem5.6}
 With the above notation, we
have
$$\|F\|_{M_d(I,\ldots, I)}\le \|f\|_{M_{d}({\bb F}_n)}
\quad{\rm and}\quad 
  \|f\|_{M_{d-1}({\bb F}_n)} \le C(d) \sup_\pi
\|F\|_{M(\pi)}$$
 where the supremum runs over all
(nontrivial)
 partitions of $[1,\ldots, d]$ into 
$K$ disjoint intervals (= blocks), with $K\le d-1$
and where $C(d)$ is a constant depending only on
$d$.
\end{lem}

\begin{proof}
The inequality $\|F\|_{M_d(I,\ldots, I)}\le
\|f\|_{M_d({\bb F}_n)}$ is essentially 
obvious by going back to the definitions, so we will now
concentrate on the 
converse direction.

\n Consider $t_1,\ldots, t_{d-1}$ in $G = {\bb F}_n$ such
that their product 
$t_1t_2\ldots, t_{d-1}$ belongs to $W_d$ \ie
$t_1t_2\ldots t_{d-1}$ can be  written as a
reduced word of the form $g_{i_1}g_{i_2}\ldots
g_{i_d}$ where 
$\{g_i\mid i\in I\}$ denotes the free generators of ${\bb F}_n$.
Since the 
letters $g_{i_1},\ldots, g_{i_d}$ remain after successive
reductions in the 
product $t_1t_2\ldots t_{d-1}$ it is easy to check that each
$t_i$ contributes a 
block of $p_i$ letters in $x$ with $\sum\limits^{d-1}_1 p_i=d$
(we allow 
$p_i=0$).  This means that when $p_i>0$, $t_i$ can be written as
a reduced word 
$x_ia_iy^{-1}_i$ with $|a_i| = p_i$ and 
when $p_i=0$  we set $a_i=e$, so 
that
\begin{equation}\label{eq5.7}
t_1t_2\ldots  t_{d-1} = a_1a_2 \ldots a_{d-1}. 
\end{equation}
Thus to each $t = (t_1,\ldots, t_{d-1})$
 as above we can associate $p(t)  = 
(p_1,\ldots, p_{d-1})$. Actually, we have
a problem here:\ this $p(t)$ is
unambiguously defined
when $d=3$ (hence we only have to
consider products of two elements). But
when $d>3$ (and thus $d-1>2$) there might be several
reductions of $t_1t_2\ldots  t_{d-1}$
leading
to the same element of $W_d$, thus there
might be several possibilities
for $p(t)$. For instance, when $d-1=3$, denoting the generators
by
$a,b,c$, the
product $abcd= (ab)(b^{-1}a^{-1}c)(c^{-1}abcd)$
(we mean here $t_1=ab $, $t_2=b^{-1}a^{-1}c$,
$t_3=c^{-1}abcd$) 
  allows
$p(t)=(0,0,4)$ but also
$p(t)=(2,0,2)$ or $p(t)=(1,0,3)$.
We prefer to ignore this difficulty for the moment
while still treating the general case, so let us
  assume $d=3$ so that $p(t)$ is always well
defined.  Moreover, if we delete the indices for
which
$p_i=0$ (and $a_i=e$), we  obtain a partition $\pi(t)$ into $k$
blocks
$(\alpha_1,\ldots, \alpha_k)$ with 
$k\le d-1$.
Then we can rewrite \eqref{eq5.7} as
\begin{equation}\label{eq5.8}
t_1t_2\ldots t_{d-1} = b_1b_2 \ldots b_k
\end{equation}
with $b_m \in W_{|\alpha_m|}$. Here we implicitly mean that the
non-reduced 
product $t_1t_2\ldots t_d$ can be viewed (just by adding
parenthesis) as a 
product
$$c_1b_1c_2b_2c_3\ldots b_k c_{k+1}$$
where each of  the intermediate products
$c_1,\ldots, c_{k+1}$ reduces to $e$. Moreover,the $k$-tuple
$(b_1,  \ldots ,b_k)$  determines a   $k$-tuple $(\hat b_1, 
\ldots ,\hat b_k)$
with $\hat b_1\in I^{\alpha_1}, \hat b_2\in I^{\alpha_2},\ldots,
\ldots,\hat b_k\in I^{\alpha_k}.$
Now fix 
$\varepsilon>0$. For any partition $\pi =
(\alpha_1,\alpha_2,\ldots, \alpha_k)$, 
we can ``factorize'' $F$ as follows:
\begin{equation}\label{eq5.9}
F(j_1,\ldots, j_k) = \eta^{\pi}_1(j_1) \ldots
\eta^{\pi}_k(j_k)\qquad (j_m \in 
I^{\alpha_m})
\end{equation}
where $\eta^{\pi}_m$ are $B(H_m,H_{m-1})$-valued functions $(H_k
= H_0 = {\bb C})$ 
such that
\begin{equation}\label{eq5.10}
\prod_m \sup \|\eta^{\pi}_m(j_m)\| \le \|F\|_{M(\pi)}
(1+\varepsilon).
\end{equation}
Let $\tilde
f$ be the function defined on $G^{d-1}$ by 
$\tilde f(t_1,\ldots ,t_{d-1})=f(t_1t_2\ldots t_{d-1})$.
Consider now the disjoint decomposition $\tilde f
=
\sum\limits_p \tilde  f_p$ where $\tilde  f_p  = \tilde f 
\cdot 1_{\{t\colon \ p(t)  = p\}}$ where the
first sum runs over all  choices of $p =
(p_1,\ldots, p_{d-1})$ with
$p_i\ge 0$ and $\sum p_i = d$.

We claim that $\|\tilde  f_p\|_{M_{d-1}(G,..,G)} \le
\|F\|_{M(\pi)}$ where $\pi$ is the partition 
associated to $p = (p_1,\ldots, p_{d-1})$ after 
removal of the empty blocks.

To prove this claim, we will produce a factorization formula for 
$\tilde  f_p(t_1,t_2,\ldots ,t_{d-1})$, namely we will
show
$$\tilde  f_p(t_1,t_2,\ldots ,t_{d-1}) = \langle\xi^p_1(t_1)
\ldots \xi^p_{d-1}(t_{d-1}) 
\delta_e, \delta_e\rangle.$$
To define $\xi^p_i(t_i)$ we must distinguish whether $p_i = 0$ or
not.

If $p_i = 0$ we set $\xi^p_i(\theta) = \lambda(\theta)\otimes  
1$ (here $\lambda(\theta)$ 
denotes as usual left translation by $\theta$ on $\ell_2(G)$). On
the other hand, if $p_i>0$ so 
that $i$ corresponds to a block $\alpha_m$ of $\pi$ with
$|\alpha_m| = p_i$, we 
write 
$$\xi^p_i(\theta) = \sum e_{x,y} \otimes \eta^{\pi}_m (\hat a)$$
where the sum runs over all ways to decompose $\theta$ as $x\cdot
a\cdot y^{-1}$ as a 
{\it reduced\/} product, with $a$ a product of generators such
that   $|a| = p_i$ and where $\hat a$ denotes the element
of $I^{\alpha_m}$ corresponding to $a$ ($x$ and $y^{-1}$ being
initial and final 
segments in the reduced word $\theta$; we allow here $x=e$ or
$y=e$). In case $\theta$ does
not admit any such decomposition (i.e. $\theta$ does not admit
any subword in $W_{p_i}$), we set  $\xi^p_i(\theta) = 0$.

\n Note that we have $\| \xi^p_i(\theta) \|\le \sup_a
\|\eta^{\pi}_m(\hat a)\|$. 
Indeed, when $\theta$ is fixed, in the various ways to write
$\theta= x\cdot a\cdot y^{-1}$ as a 
{\it reduced\/} product as above, all the $x$'s appearing will be
distinct since they have different length,
and similarly all the $y$'s will be distinct, so the various
operators $ e_{x,y} \otimes \eta^{\pi}_m (\hat a)$
have both orthogonal ranges and orthogonal domains, so that the
norm of their sum
is majorized by the maximum norm of each term.

\n  A (tedious but) straightforward verification shows that if
$t_1t_2\ldots t_{d-1} 
\in W_d$ with $p(t) = p$, and if $t_1t_2\ldots t_{d-1} = b_1b_2
\ldots b_k$ as described in \eqref{eq5.8}, then we have using \eqref{eq5.9}
\begin{align}\label{eq5.11}
\langle \xi^p_1(t_1)\ldots \xi^p_{d-1}(t_{d-1})
\delta_e, 
\delta_e\rangle &= \eta^{\pi}_1(\hat b_1)\ldots 
\eta^{\pi}_k(\hat b_k)\\
&= F(\hat b_1,\ldots,\hat b_k)\nonumber\\
\text{whence by \eqref{eq5.6} and \eqref{eq5.8}}\qquad &= f(t_1\ldots
t_{d-1}).\nonumber
\end{align}
Moreover, if $p(t) \ne p$ the left side of \eqref{eq5.8} vanishes.

\n Indeed, if that left side
is non zero, then we must have 
$$t_1\cdots t_{d-1}= [x_1 .a_1.y_1^{-1}][ x_2.a_2.y_2^{-1}]
\cdots [x_{d-1} .a_{d-1}.y^{-1}_{d-1}]
=a_1a_2 \cdots a_{d-1}
$$  with
$a_i\in W_{p_i}$ if $p_i>0$,  and
$a_i=e$ otherwise, $a_i$ being a subword of $t_i$,
 in such a way that the product of all the terms
figuring in between two successive $a_i$'s with $p_i>0$ reduces
to $e$,
as well as the product of all the terms preceding the first $a_i$
with $p_i>0$,
and that of all the terms after the last $a_i$ with $p_i>0$.
 This implies $p(t) = p$ and deleting the $a_i$'s equal to $e$ we
obtain
 $t_1\cdots t_{d-1}=b_1\cdots b_k$ and \eqref{eq5.8} is then easy to
check.

\n Thus we have the 
announced factorization of $\tilde  f_p$; the latter implies 
$$\|\tilde  f_p\|_{M_{d-1}(G,..,G)}
 \le \prod \sup\|\eta^{\pi}_m\| \le \|F\|_{M(\pi)}
(1+\varepsilon).$$ 
Using $ \tilde f = \sum
\tilde  f_p$, this yields  
$\|  f\|_{M_{d-1}(G)}=\|\tilde
f\|_{M_{d-1}(G,..,G)} \le C_d \|F\|_{M(\pi)}\allowbreak 
(1+\varepsilon) $ 
(here $C_d $ is the number of possible $p$'s) thus completing the
proof of
the  lemma, at least in the case $d=3$.
Since there are only four possibilities for
$p$ (namely $(3.0),(0.3), (1,2)$,  and $(2.1)$) we obtain  
$C_3\le 4 $.

 Now in the general case, the problem is that,
for each $t=(t_1,\ldots, t_{d-1})\in G^{d-1}$
such that $t_1\ldots t_{d-1}\in W_d$, 
there is a multiplicity of possible $p(t)'s$
(or   of possible associated partitions $\pi(t)$):
each such $t$  admits
$N(t)$ possible distinct
$p(t)'s$. However,  we  of course
have a bound for this: $1\le N(t)\le N_d$ where
the upper bound
$N_d$ depends only on $d$. If  
$t_1\ldots t_{d-1}\not\in W_d$, we set $N(t)=0$. Then,
we think of
$p(t)$ as a multivalued function and we define
$$\tilde  f_p  = \tilde f 
\cdot 1_{\{t\colon \ p\in p(t)    \}}.$$
Then the preceding shows again that 
\begin{equation}\label{eq5.12}
 \|\sum \tilde  f_p\|_{M_{d-1}(G,..,G)} \le
 C_d \|F\|_{M(\pi)}
(1+\varepsilon) ,
\end{equation}
but, since the sum is no longer
disjoint, we have 
\begin{equation}\label{eq5.13}
 {\forall t\in G^{d-1}} \quad  \sum_p \tilde  f_p(t) =  N(t)
\tilde f(t). 
\end{equation}
Consider now the special case when $F$ is identically equal to 1. 
Note that $\|F\|_{M(\pi)} \le 1$ 
and $ N(t)
\tilde f(t)  =N(t)$ in this case. Thus, the preceding identity
and \eqref{eq5.12} shows
that the function
$N\colon \ G^{d-1}\longrightarrow {\bb R}$ is in
$M_{d-1}(G,\ldots, G)$ with norm $\le C_d$.
To conclude, we will mupliply \eqref{eq5.13} by a 
function equal to $1/N$ on the support of $N$ and we will bound
its norm in $M_{d-1}(G,\ldots,
G)$ by a constant $C_d'$. Since $M_{d-1}(G,\ldots,
G)$ is a Banach algebra for the pointwise product, this will
yield the desired result. (Alternately, we could use
a disjointification trick, as above for   \eqref{eq2.151}.)

Let $P$ be a polynomial  such that $P(k) = \frac1k$ for all
$k=1,2,\ldots, N_d$. To fix ideas, we let $P$ be determined 
by Lagrange interpolation. Since
 $M_{d-1}(G,\ldots,
G)$ is a Banach algebra,  $P(N) \in M_{d-1}(G,\ldots, G)$ and
since $P$ depends only on $d$,  we have
$\|P(N)\|_{M_{d-1}(G,\ldots, G)} \le C'_d$ for some $C'_d$
depending only on
$d$. Then we can write
$$P(N) \cdot \sum_p\tilde f_p = P(N)\cdot N\tilde f
= \tilde f$$ 
hence we conclude
$$\|f\|_{M_{d-1}(G)} = \|\tilde f\|_{M_{d-1}(G,\ldots, G)}\le \|
P(N) \|_{M_{d-1}(G,\ldots, G)}
\|\sum \tilde  f_p\|_{M_{d-1}(G,..,G)}$$
$$\le C'_dC_d
\|F\|_{M(\pi)}
(1+\varepsilon) .\qquad\qed$$
 \renewcommand{\qed}{}\end{proof}

\begin{proof}[Proof of Theorem \ref{thm5.1}]
   Assume $M_d({\bb F}_\infty) = 
M_{d-1}({\bb F}_\infty)$. Then there must exist a constant $C'$
such that for all $f$  in $M_{d-1}({\bb F}_\infty)$ we have
$$\|f\|_{M_d({\bb F}_\infty)} \le C'\|f\|_{M_{d-1}({\bb F}_\infty)}.$$
Then by Lemma \ref{lem5.6} we find that $\Phi^{-1}\colon \
\bigcap\limits_\pi M(\pi)\to  M_d(I,\ldots, I)$
is bounded, which contradicts Lemma~\ref{lem5.3} for any
$d>1$.

\n  Now,
let  $F_{d,n}$ be as   in the remark following
Lemma  \ref{lem5.5}, and let $f_{d,n}$ be 
defined  by
$$f_{d,n}(g_{i_1}g_{i_2}\ldots g_{i_d}) = F_{d,n}(i_1,\ldots,
i_d)$$
and $f_{d,n}(t) = 0$ if
 $t\notin W_{d,n}$. Then by the latter remark
  and by Lemma \ref{lem5.6} we have 
$$n^{\frac{d-1}{2}} \le \|f_{d,n}\|_{M_d({\bb F}_\infty)},$$
and also by Lemma \ref{lem5.6} and Lemma \ref{lem5.4}
$$\|f_{d,n}\|_{M_d({\bb F}_\infty)} \le C(d)
\sup_{\pi}\|F_{d,n}\|_{M(\pi)} \le C(d)n^{\frac{d-2}{2}}. $$
This complete the proof. 
\end{proof}

Note that in the special case $d=3$, we   obtain a very explicit
example:
Namely the function $f_{3,n}$ supported on $W_{3,n}$ and defined
there for $1\le p,q,r\le n$ by
$$f_{3,n}(g_p g_q g_r)= \exp(i(p+r)q/n)  ,$$
 satisfies
$$ n\le \|f_{3,n}\|_{M_3(G)}\quad \text{but} \quad
\|f_{3,n}\|_{M_2(G)} \le 4 n^{1/2}.$$

\begin{rem}\label{rem5.7} 
 The proof of Lemma \ref{lem5.6} can be modified to
show that 
$\|f\|_{M_d({\bb F}_n)} \le C'(d) \|F\|_{M_d(I,\ldots, I)}$ for some
constant $C'(d)$ 
depending only on $d$. In particular, we have
$\|1_{W_{d,n}}\|_{M_d({\bb F}_n)}\le 
C'(d)$ and hence $\|1_{W_d}\|_{M_d({\bb F}_\infty)}\le C'(d)$.
\end{rem}

Note however:

\begin{pro}\label{pro5.8}
  Any function $f\colon \ {\bb F}_\infty \to 
{\bb C}$, supported on $W_d$, that is in $M_d({\bb F}_\infty)$ must
necessarily be in 
$UB({\bb F}_\infty)$. 
\end{pro}
 
\begin{proof}
Indeed, $f_{|W_d}\colon \ W_d\to {\bb C}$ admits an extension to
a 
function $\hat f\colon \ {\bb F}_\infty\to{\bb C}$ that is in
$C^*({\bb F}_\infty)^* = 
B({\bb F}_\infty)$. This follows from \cite[Corollary 8.13]{P8}. 
Now, since, by Remark \ref{rem5.7}, $1_{W_d}$ belongs to $UB({\bb F}_\infty)$,
the pointwise product 
$f=\hat f \cdot 1_{W_d}$ also belongs to $UB({\bb F}_\infty)$.
\end{proof}

This shows in particular that   a function in 
$M_2({\bb F}_\infty)\backslash M_3({\bb F}_\infty)$ cannot be supported
on $W_2$.
Thus the above example $f_{3,n}$ supported on $W_3$ appears
somewhat ``minimal".

Let $G$ be any free group. Recall that we denote by ${\cl W}(d)$
the set of all 
words of length $d$ in the generators and their inverses. Note
that the inclusion 
$W_d \subset  
{\cl W}(d)$ is strict. We chose to concentrate on $W_d$ (rather
than on ${\cl W}(d)$) 
because then the ideas are a bit simpler and Lemma \ref{lem5.6} is
somewhat prettier in 
that case:\ Indeed, that lemma identifies the spaces $\{f\in
M_d(G)\mid 
\text{supp}(f)\subset W_d\}$ and $\{f\in M_{d-1}(G)\mid
\text{supp}(f)\subset 
W_d\}$ with two distinct spaces of functions on $G^d$, thus
reducing, in some 
sense, a problem in harmonic analysis to one in functional
analysis. However, 
most of our results hold with suitable modification for functions
with support in 
${\cl W}(d)$. We will merely describe them with mere indication
of proof.

 Fix $d\ge 
1$ and let $k\le d$. Let $f\colon \ G\to {\bb C}$ be a function
supported on ${\cl 
W}(d)$. Let $\pi$ be a partition of $[1,\ldots, d]$ in $k$
disjoint consecutive 
blocks (intervals) $\alpha_1,\ldots,\alpha_k$ so that $|\alpha_1|
+\cdots+ 
|\alpha_k| = d$ and $|\alpha_i|\ge 1$ for all $i=1,\ldots, k$. We
will denote by 
$k(\pi)$ the number of blocks, i.e.\ we set $k(\pi) = k$. We
define 
$f_\pi\colon \ {\cl W}(|\alpha_1|) \times\cdots\times {\cl
W}(|\alpha_k|)\to 
{\bb C}$ by
$$f_\pi(x_1,x_2,\ldots, x_k) = f(x_1x_2\ldots x_k).$$
Note that $f_\pi(x_1,\ldots, x_k)=0$ if the product $x_1x_2\ldots
x_k$ is not a 
reduced word, since then it has length $<d$. For any function
$F\colon \ {\cl 
W}(|\alpha_1|)\times\cdots\times {\cl W}(|\alpha_k|) \to {\bb C}$.
we denote again
\begin{equation}\label{eq5.14}
\|F\|_{{\cl M}(\pi)} = \|F\|_{M({\cl W}(|\alpha_1|),\ldots,
{\cl 
W}(|\alpha_k|))}.
\end{equation}
The preceding proofs (mainly Lemma \ref{lem5.6}) then yield

\begin{thm}\label{thm5.9}
 With the preceding notation, we have for
any function $f$ 
with support in ${\cl W}(d)$ and for any integer $K\geq 1$
$$\sup_{k(\pi)\le K} \|f_\pi\|_{{\cl M}(\pi)} \le \|f\|_{M_K(G)}
\le C(d,K) 
\sup_{k(\pi)\le K} \|f_\pi\|_{{\cl M}(\pi)}$$
where $C(d,K)$ is a constant depending only on $d$ and $K$.
\end{thm}

\begin{rk}
 Let $\pi_0$ be the partition of $[1,\ldots, d]$
into 
singletons, so that $k(\pi_0) = d$. When $K\ge d$, we have for
any $f$ as in 
\eqref{eq5.14}
\begin{equation}\label{eq5.15}
\sup_{k(\pi)\le K} \|f_\pi\|_{{\cl M}(\pi)} =
\|f_{\pi_0}\|_{{\cl 
M}(\pi_0)}. 
\end{equation}
Indeed, we have $k(\pi_0)\le K$ and moreover it is easy to see
using \eqref{eq5.3} that 
if a partition $\pi$ is less fine than another one $\pi'$ (i.e.\
every block in 
$\pi$ is the union of certain blocks of $\pi'$) we have
$$\|f_\pi\|_{{\cl M}(\pi)} \le \|f_{\pi'}\|_{{\cl M}(\pi')}.$$
Since any $\pi$ is less fine that $\pi_0$, \eqref{eq5.15} follows
immediately.
\end{rk}

Thus, for all $K\ge d$, the norms $\|f\|_{M_K(G)}$ are equivalent
on functions 
$f$ with support in ${\cl W}(d)$. Indeed, by \eqref{eq5.14} they are
equivalent to 
$\|f_{\pi_0}\|_{{\cl M}(\pi)}$. In sharp contrast when $K<d$, in
particular 
when $K=d-1$, they are no longer equivalent.

\begin{cor}\label{cor5.10}
Let $G$ be any group containing a (non-Abelian)
 free subgroup. Then, 
for any $\theta>1$ and any $c>\theta$ we have
$B_c(G)\not=B_\theta(G)$,
consequently there is a representation $\pi\colon \ 
G\to B(H)$ with $|\pi|\le c$ that is not similar to any
representation $\pi'$ with 
$|\pi'|\le\theta$.
\end{cor}

\begin{proof}
Assume by contradiction that the conclusion fails. Then, a
fortiori it must also 
fail for $G={\bb F}_\infty$, by an ``induction'' argument as in
Proposition \ref{pro0.5}. 
Hence, by Theorem \ref{thm2.9}, we obtain $M_d({\bb F}_\infty) =
M_{d+1}({\bb F}_\infty)$ for some 
$d$, contradicting Theorem \ref{thm5.1}.
\end{proof}

\begin{rk}
 In the case of $G = {\rm SL}_2({\bb R})$, Michael
Cowling showed me a 
very concrete proof (that he attributed to Haagerup) of the
conclusion of 
Corollary \ref{cor5.10}. That proof uses the estimates of Kunze and Stein
from \cite{KS1}, the 
Bruhat decomposition of ${\rm SL}_2({\bb R})$ and the amenability of
the subgroup of 
triangular matrices in ${\rm SL}_2({\bb R})$.
\end{rk}

\begin{rk}
 In \cite{P6}, we study the same
question as in this section
but for $G={\bb N}$. We answer a related question of
Peller concerning power bounded operators
(=uniformly bounded representations of
$G={\bb N}$), by showing
$$M_2({\bb N})\not= M_{3}({\bb N}).$$
On the other hand, the main result of 
\cite{KLM} implies that if $H=\ell_2$ we have for any
$d$
$$M_d({\bb N};H)\not= M_{d+1}({\bb N};H).$$
However, the same question for $H={\bb C}$ remains
open when $G={\bb N}$.
\end{rk}

\n{\it Acknowledgement.} 
This is a much  expanded version of a manuscript 
initially based on the notes  of my talk
at the 1998 Ponza conference for A. Fig\`a-Talamanca's 60-th birthday.
I am very grateful to
Michael Cowling for useful conversations and for
providing Remark \ref{rem0.72}.
I am indebted to Marius Junge for Lemma \ref{lem5.5}.
 I also  thank Slava
Grigorchuk for useful remarks on the text.

\end{document}